\documentclass[11pt, final]{amsart}
\usepackage[english]{babel}
\usepackage{latexsym}
\usepackage{amssymb}
\usepackage{amscd}
\usepackage[cp850]{inputenc}
\usepackage[mathscr]{eucal}

\usepackage[all]{xy}

\newcommand*{\xdash}{\rule[0.5ex]{0,8em}{0.5pt}\vert\,}
\newcommand{\adef}{\begin{defin}}
\newcommand{\zdef}{\end{defin}}
\newcommand{\CAT}{\mathbf}

\theoremstyle{plain}
\newtheorem{theorem}{Theorem}[section]

\newtheorem{defin}[theorem]{Definition}
\newtheorem{proposition}[theorem]{Proposition}
\newtheorem{cor}[theorem]{Corollary}
\newtheorem{lemma}[theorem]{Lemma}

\theoremstyle{remark}

\newcommand{\R}{\mathbb{R}}

\newcommand{\coker}{\mathrm{coker}}

\def\Ext{\operatorname{Ext}}

\def\co{\operatorname{co}}

\def\Hom{\operatorname{Hom}}
\newcommand{\aproof}{\begin{proof}}
\newcommand{\zproof}{\end{proof}}

\newcommand{\To}{\longrightarrow }

\newcommand\PB{{\mathrm{PB}}}

\begin{document}
\title[Categorical Banach space theory II]{The hitchhiker guide to Categorical\\ Banach space theory. Part II.}

\author{Jesus M. F. Castillo}
\address{Universidad de Extremadura, Instituto de Matem\'aticas Imuex, 
Avenida de Elvas s/n, 06011 Badajoz, Spain} \email{castillo@unex.es}

\thanks{2020 AMS classification: 46M15, 46M18, 18-02, 18A40}
\thanks{This research was supported in part by MINCIN Project PID2019-103961GB-C21 and project IB20038 de la Junta de Extremadura.}

 \maketitle

\begin{abstract} What has category theory to offer to Banach
spacers? In this second part survey-like paper we will focus on very much needed advanced categorical and homological elements, such as Kan extensions, derived category and derived functor or Abelian hearts of Banach spaces.\end{abstract}

\section*{HMBST. So it goes.}

\emph{The Hitchhiker Guide to Categorical Banach Space Theory. Part I} [Extracta Math. 25 (2010) 103--149], from now on referred to as \texttt{HHI} was written to display the basic elements from categorical algebra transplanted to Banach space theory:
\begin{itemize}
\item The definition of category, functor and natural
transformation.
\item Pullbacks and pushout construction.
\item Limits and colimits.
\item Adjoint functors.
\item Duality
\end{itemize}
As it was promised there ``the study of derived functors conforms what is called homology theory, and will be treated in Part II". Ten years later, it is showtime! In the span between \texttt{HHI} and the present moment, quite a few things happened: one of them is the writing and imminent publication of the book [F. Cabello S\'anchez, J.M.F. Castillo, \emph{Homological Methods in Banach Space Theory}, Cambridge Studies in Advanced Mathematics, Cambridge Univ. Press] \cite{hmbst}, from now on referred to as HMBST. This book contains, as the title clearly states, the lion's share of what (the authors know that) is currently known about homological methods in Banach space theory. More precisely, the following topics, always with the focus set on Banach and quasi Banach spaces, are treated there: The definition of category, functor and natural transformation, pullbacks and pushouts constructions, exact sequences, the functor $\Ext$, derivation of functors and why $\Ext$ is the derived functor of $\mathfrak L$ in the category of $p$-Banach spaces, the construction of the long homology sequence up to $\Ext$ terms, limits and colimits, adjoint functors, natural equivalences of $\Ext$ and Fra\"iss\'e limits. Since it is absurd to ignore that, but also because we will not be able to provide a better exposition than the one presented in HMBST, this survey continues where HMBST stopped. We have also made every possible effort not to continue where \cite{pfiff} stopped.\\

A point HMBST attempted to make clear is that the study exact sequences of Banach spaces needs ---actually, necessarily, peremptorily needs--- dealing with exact sequences of quasi Banach spaces. There is no escape to that. Thus, if one accepts that homology deals with exact sequences,
whatever attempt to study homology in the category of Banach spaces requires coming to terms with the idea of making homology in the category of quasi Banach spaces. But not all things that can be done in Banach spaces can be done in quasi Banach spaces. And this is, categorically speaking, a catastrophe. Throughout this paper it is our plan to display several categorical tools that could help to overcome, surround or circumvent such difficulties.

Categories will be labeled $\CAT{A,B,C,D, E, F}\dots$, and functors $\mathcal{A,B,C,D, E, F}\dots$. Those of Banach and quasi Banach spaces (and operators) will be called $\mathbf{Ban}$ and $\mathbf{QBan}$; and that of vector spaces and linear maps that will be called $\mathbf{Vect}$. The categories of Banach (resp. quasi Banach) spaces and contractive operators will be called $\mathbf {Ban_1}$ (resp. $\mathbf{QBan_1}$).

\section{The categories $\CAT C^{\CAT D}$}\label{catCD}

Given a category $\CAT C$ and a small category (say, one that can be understood as a diagram; or, one whose objects and morphisms are sets) $\CAT D$, the product category $\CAT C^{\CAT D}$ has as objects
the functors $\CAT D \to \CAT C$ and the corresponding natural
transformations as morphisms. There is an obvious
covariant diagonal functor $\Delta: \CAT C \to \CAT C^{\CAT D}$: the
image of an object is the functor $\CAT D \to \CAT C$ sending all
objects to $X$ and all arrows to the identity of $X$. Which is the image of an arrow is clear. This uplifting of the action from $\CAT C$ to $\CAT C^{\CAT D}$ can be enriching. Think about $\mathbf{Ban}$. A conceptually valid reason to jump  from $\mathbf{Ban}$ to $
{\mathbf{Ban}}^{\mathbf{Ban}}$ is to provide a solid support to the Eilenberg-MacLane program (which, as explained in HHI, means that functors and natural transformations is all there is in a mathematical theory). Thus, Banach spaces must be interpreted as functors and operators would then be natural transformations. Let us attempt identifying the category \textbf{Ban} of [Banach spaces +
operators] as a full subcategory of the category ${\mathbf{Ban}}^{\mathbf{Ban}}$ of [Banach functors +  Banach natural
transformations] (see below for unexplained notation). Recall that a \emph{Banach functor} is a functor $\mathcal F: \textbf{Ban} \to
\textbf{Ban}$. It is often required that $\mathcal F$ is linear,
with the meaning that $\mathcal F (\lambda T+S)= \lambda \mathcal
F (T) + \mathcal F (S)$; also, it will also be sometimes required
that it is norm decreasing, with the meaning $\|\mathcal F (T) \|
\leq \|T\|$. A  number of reasonably interesting Banach functors were presented in HHI:
\begin{itemize}
\item Given a Banach space $Y$ the contravariant  $\mathfrak L^Y$
functor defined by $\mathfrak  L^Y(X)= \mathfrak L(X,Y)$ and
$\mathfrak  L^Y(T)(S)= ST$. The choice $Y=\R$ gives the duality functor. \item Given a Banach
space $X$ the covariant $\mathfrak  L_X$ functor defined by
$\mathfrak  L_X(Y) = \mathfrak L(X,Y)$ and $\mathfrak  L_X(T)(S) = TS$. The choice $X=\R$ gives the identity. \item Given a
Banach space $X$ the covariant $\otimes_X$ functor defined by
$\otimes_X(Y)= X \widehat{\otimes}_\pi Y$ and $\otimes_X(T) = 1_X
\otimes T$. \item Semadeni's covariant Banach-Mazur
functor \cite{semaprace} that assigns to a Banach
space $X$ the space $C(B_{X^*})$; and to a contractive operator $T: X
\To Y$ the operator $f \to fT^*$.
\item The ultraproduct functor $X\to X_{\mathscr U}$.
\item The remainder functor $X\to X^{**}/X$.
\item The covariant functors assigning to a Banach space $X$ the space $\ell_p(X)$ of
$p$-summable sequences or the space $c_0(X)$ of norm null sequences and the naturally induced operators.
\item The Grothendieck-Pietsch functors that
assign to a Banach space $X$ the space $\ell_p^w(X)$ of weakly
$p$-summmable sequences on $X$ or the space $c_0^w(X)$ of weakly null sequences.
\end{itemize}

In addition to that, there are a couple of distinguished functors between $\mathbf{Ban}$ and $\mathbf{QBan}$.

\begin{itemize}
\item The forgetful functor $\square: \mathbf{Ban}\To \mathbf{QBan}$ that simply forgets the Banach structure and just leaves the quasi Banach structure.
\item The Banach envelope functor $co: \mathbf Q\To \mathbf B$ that associates to each quasi Banach space $Q$ with unit ball $B_Q$ its Banach envelope $co(Q)$, which is the completion of $Q$ endowed with the norm whose unit ball is the convex envelope of $B_Q$ (and, if the dual of $Q$ does not separate points, making the necessary quotient). This is a functor \cite[1.1]{hmbst}: the canonical map $\imath_Q: Q\To co(Q)$ has the property that whenever $Y$ is a Banach space every operator $Q\To Y$ factors through $\imath_Q$.
    \end{itemize}

A natural transformation $\eta: \mathcal F\To \mathcal G$ is called a  Banach natural transformation if
$\|\eta\| = \sup_X \|\eta_X: \mathcal FX\To \mathcal GX\|<\infty$. The Banach space of Banach natural transformations between two Banach functors $\mathcal F, \mathcal G$ will be denoted $[\mathcal F, \mathcal G]$.\\

There are however many ways to reinterpret $\mathbf{Ban}$ inside $\mathbf{ Ban}^{\mathbf{ Ban}}$. In a sense, everything one has to do is to replace Banach spaces by Banach functors, which can be done in many ways, and operators by Banach
natural transformations. More precisely, one has to define an interpretation functor $\delta: \mathbf{Ban}\To \mathbf{ Ban}^{\mathbf{ Ban}}$.
For instance (see HHI for details):
\begin{itemize}
\item If $\delta(X ) = \otimes_X$ then $[\otimes_A, \otimes_B] = \mathfrak L(A,B)$; which shows this is a faithful representation.
\item If $\delta(X ) = \mathfrak L_X$ then $[\mathfrak L_A, \mathfrak L_B] =
\mathfrak L(B,A)$, so this is also a faithful representation.
\item If $\delta(X ) = \mathfrak L^X$ then $[\mathfrak L^A, \mathfrak L^B] =
\mathfrak L(A,B)$, so this is also a faithful representation.
\item Of course, one could set $\delta = \Delta$, in which case one also has $\mathfrak L(X,Y) = [\Delta(X),\Delta(Y)]$.
\end{itemize}

In this way, when working on $\mathbf{Ban}$ the action does actually takes place in ${\mathbf{ Ban}}^{\mathbf{ Ban}}$, and one must thus be ready to transplant our notions to ${\mathbf{ Ban}}^{\mathbf{ Ban}}$. How to do that? An idea is to trust in Kan extensions. Naively said, once an interpretation $\delta$ has been established, a Kan extension of a Banach functor $\mathcal F: \mathbf{Ban} \to
\mathbf{Ban}$ is a functor $\mathcal F^K:
\mathbf{Ban}^{\mathbf{Ban}} \longrightarrow
\mathbf{Ban}^{\mathbf{Ban}}$ such that $\mathcal F^K \delta =
\delta \mathcal F$. To understand what is going on here, we need to recall a couple of notions.

\section{Adjoint functors}

The adjointness notion, also due to Kan \cite{kan}, was explained in Part I and will be explained now. Given two covariant functors $\mathcal F:
\CAT{A} \to \CAT{B}$ and $\mathcal G: \CAT{B} \to \CAT{A}$ one says that $\mathcal F$ is a left adjoint of $\mathcal G$, written
$\mathcal F \xdash \mathcal G$ (and, consequently, $\mathcal G$ is a right adjoint of $\mathcal F$) when for every objects $A$
of $\CAT{A}$ and $B$ of $\CAT{B}$ there is a natural equivalence
$$\eta: \mathrm{Hom}_\mathbf{B}(\mathcal F(\cdot), \cdot)  \longrightarrow \mathrm{Hom}_\mathbf{A}(\cdot, \mathcal G(\cdot))$$
between the functors $A,B \to \mathrm{Hom}_\mathbf{B}(\mathcal F(A), B)$ and $A,B \to \mathrm{Hom}_\mathbf{A}(A, \mathcal G(B))$. As Semadeni mentions \cite[p.217]{semabook}, sometimes during proofs it is more convenient to use an equivalent formulation of the adjunction relation: it turns out that $\mathcal F \xdash \mathcal G$ if and only if there are natural transformations $\varepsilon: 1 \to \mathcal G \mathcal F$ and $\delta: \mathcal F \mathcal G\to 1$ such that
$\delta \mathcal F\circ \mathcal G\varepsilon=1$ and $\mathcal G\delta \circ \varepsilon \mathcal G=1$ as natural transformations
$$\xymatrix{\mathcal F \ar[r]^{\mathcal F\varepsilon}&\mathcal F\mathcal G\mathcal F\ar[r]^{\delta \mathcal F}& \mathcal F}\quad\quad\xymatrix{\mathcal F \ar[r]^{\varepsilon \mathcal G}&\mathcal G\mathcal F \mathcal G\ar[r]^{\mathcal G\delta }& \mathcal G}$$

In other words, if we decide to ignore the natural transformations involved and replace them by equalities, the idea behind this notion is that
$\mathcal F \xdash \mathcal G$ when $\mathrm{Hom}_\mathbf{B}(\mathcal F(A), B) =\mathrm{Hom}_\mathbf{A}(A, \mathcal G(B))$. A second interpretation is only possible with care: one is certainly saying that $\mathcal G$ and $\mathcal F$ are, somehow, inverses, namely $1= \mathcal G \mathcal F$ and $\mathcal F\mathcal G =1$ in the sense that there are natural transformations $1\to \mathcal G \mathcal F$ and $\mathcal F\mathcal G \to 1$ (but not  $\mathcal G \mathcal F=1$ or $1=\mathcal F \mathcal G$ !).

Just to give an example at hand, let us recall from \cite[4.6.1 (b)]{hmbst} that $co \;\xdash \square$.

\section{Limits} This topic is covered both in \texttt{HHI} and HMBST. Briefly said, given a functor $\mathcal F\in \CAT C^{\CAT D}$, its limits
$\lim \limits_\rightarrow \mathcal F$ and $\lim\limits_\leftarrow\mathcal F$, if they exist, are objects of $\CAT C$ provided with arrows $\alpha_d: \mathcal F(d)\to \lim\limits_\leftarrow\mathcal F$ (resp.
$\beta_d: \lim \limits_\rightarrow\mathcal F \to\mathcal F(d)$) making a commutative diagram, and being universal with respect to this property; namely, given any other object $c$ of $\CAT C$ provided with arrows $\alpha'_d: \mathcal F(d)\to c$ making the diagram commutative, there is a unique arrow $\alpha: \lim\limits_\rightarrow \mathcal F\to c$ making all diagrams commutative (resp. given any other object $c$ of $\CAT C$ provided with arrows $\alpha'_d: c\to \mathcal F(d)$ making the diagram commutative, there is a unique arrow $\beta: c\to \lim\limits_\leftarrow \mathcal F$ making all diagrams commutative). Enthusiastic readers will have surely observed that when limits exist then $\lim \limits_\rightarrow$ and $\lim \limits_\leftarrow$ can be understood as functors $\CAT C^{\CAT D} \To \CAT C$, which turn out to be the adjoints of $\Delta: {\CAT C} \to \CAT C^{\CAT D}$; precisely:
$$\Delta \;\; \xdash \;\;\lim\limits_\leftarrow F\quad\quad \quad \mathrm{and} \quad\quad \quad\lim_{\rightarrow}\mathcal F\;\; \xdash \;\; \Delta$$
A lively discussion about which one,  $\lim \limits_\rightarrow$ or $\lim \limits_\leftarrow$, is \emph{the} limit and which one \emph{the} co-limit is in \texttt{HHI}. Here we will call, when the occasion demands it, $\lim\limits_\rightarrow$ the inductive limit and $\lim\limits_\leftarrow$ the projective limit. A simple example can be provided just for fun:

\begin{proposition} In the category $\mathbf{Ban_1}$, every Banach space is an inductive limit of its finite dimensional subspaces.
\end{proposition}

HINT: If the reader feels that she (or he) \emph{really} needs a proof, do it. Or else, read it \cite[I, 1.20] {ciglosmich}.

\section{Comma categories}

Furnishing us with this tool is like putting a gun in a monkey's hand. Or, worse, a typewriter.
Assume that one has a category $\CAT C$ and a larger category $\CAT D$ and pick an object $D$ of $\CAT D$.

\adef The comma category $(\CAT C \Lsh D)$ has as objects pairs $(X,f)$ formed by an object $X$ of $\CAT C$ and an arrow $f: D\To X$, and whose arrows $\varphi: (X,f)\To (Y, g)$ are arrows $\varphi: X\To Y$ so that $\varphi g=f$. The comma category $(\CAT C \Rsh D)$ has as objects pairs $(X,f)$ formed by an object $X$ of $\CAT C$ and an arrow $f: X\To D$, and whose arrows $\varphi: (X,f)\To (Y, g)$ are arrows $\varphi: X\To Y$ so that $g\varphi =f$.\zdef

Any functor $\mathcal F: \CAT C \To \CAT E$ induces functors $\mathcal F^D: (\CAT C \Rsh D) \To \CAT E$ given by $\mathcal F^D(X,f)= \mathcal F(X)$ and another functor $\mathcal F_D: \CAT (\CAT C \Lsh D) \To \CAT E$ given by $\mathcal F_D(X,f)= \mathcal F(X)$. Paraphrasing Roberto Begnini, when a man with a gun meets a man with a pen, the man with the gun is done. Let us show why the monkey was more dangerous with the comma category:

\begin{proposition} Every Banach space $X$ is an inductive limit $X=\lim\limits_{\To} \ell_1^n$ in $\mathbf{Ban_1}$.
\end{proposition}
\begin{proof} Think of $\mathbf {Ban_1}$ as ``the larger category" inside which one will work. Construct the comma category $(\ell_1^\bullet \Rsh X)$ where $\ell_1^\bullet$ represents here the (full sub-) category (of $\mathbf {Ban_1}$) whose objects are finite dimensional $\ell_1^n$ spaces and let $\mathcal F: (\ell_1^\bullet, \Rsh X)\To \mathbf{Ban_1}$ the obvious forgetful functor that associates to an object $(\ell_1^n, \varphi)$ the space $\ell_1^n$. Let us show that $X=\lim\limits_{\To} \mathcal F( \ell_1^d, \varphi_d) = \lim\limits_{\To} \ell_1^d$. This amounts showing that whenever one has a Banach space $Y$ and another commuting family of arrows $\psi_d: \ell_1^d\To Y$ (i.e., when $\alpha_{d,d'}: \ell_1^d\To \ell_1^{d'}$ in the comma category ---$\varphi_{d'}\alpha_{dd'}=\varphi_d$--- as in the diagram below) then it is possible to define a contractive operator $T:X\To Y$ so that $T\varphi_d = \psi_d$. Of course that given a commuting diagram
$$
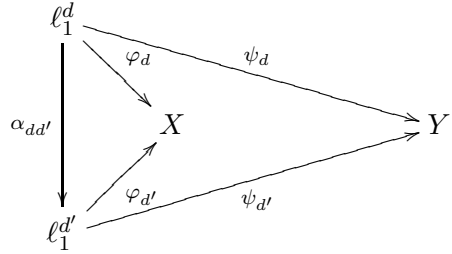
\xymatrix{
\ell_1^d \ar[dr]^{\varphi_d}\ar[dd]_{\alpha_{dd'}}\ar[drrrr]^{\psi_d}&&&\\
&X&&&Y\\
\ell_1^{d'} \ar[ur]_{\varphi_{d'}}\ar[urrrr]_{\psi_{d'}}&}$$
there is not, necessarily, a $T:X\To Y$ making the whole diagram commutative. But that is on a one-by-one basis. When the families $(\varphi_d)$ and $(\psi_d)$ involve all choices, such $T$ exist. To prove it let us show that whenever $\varphi_d(u)=\varphi_{d'}(u')$ then
$\psi_d(u)=\psi_{d'}(u')$. If so, given $x\in X$ such that $x=\varphi_d(u)$ the point $T(x)=\psi_d(u)$ is well defined.

Now, pick $x\in X$ and $\phi_d: \ell_1^d\To X$ so that $\phi_d(u)=x$. Since $\phi_d$ is contractive, $\|x\|\leq \|u\|_1$. Let $L$ be the one dimensional Banach space. When it is not confusing we simplify $\ell_1^d$ to just $d$ and define $\phi[d,\phi_d, u, x]: L\To X$ by $\phi[d, \phi_d, u, x](1)= x/\|u\|_1$. The pair $(L, \phi[d,\phi_d, u,x])$ is an object in the comma category, and $\alpha: (L, \phi[d,\phi_d, u,x]) \To (\ell_1^d, \phi_d)$ given by $\alpha(1)= u/\|u\|_1$ defines a morphism of the comma category. So we are all friends here, right? And with this we want to say that when the family $(\psi_d)$ enters the game then $\psi_d=\psi_d\alpha$ which, in particular, yields $$\psi_d(u)= \psi_d\alpha(\|u\| 1) = \|u\|\psi_d\alpha(1) = \|u\| \psi_{(L, \phi[d,\phi_d, u,x])}(1).$$

Assume that $\varphi_d(u)=\varphi_{d'}(u')$ with $\|u\|\geq \|u'\|$ (the other case is similar). The commutativity of the diagram
$$\xymatrixrowsep{0,5cm}\xymatrixcolsep{3cm}
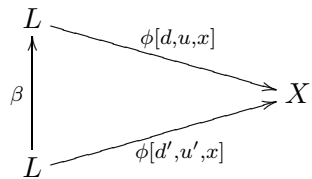
\xymatrix{
L \ar[dr]^{\phi[d,u,x]}&\\
&X\\
L \ar[uu]^\beta \ar[ur]_{\phi[d',u',x]}}$$
in which $\beta(1) = \|u\|/\|u'\|$ provides that also
$$\xymatrixrowsep{0,5cm}\xymatrixcolsep{3cm}\xymatrix{
L \ar[dr]^{\psi_{(L, \phi[d,u,x])}}&\\
&Y\\
L \ar[uu]^\beta \ar[ur]_{\psi_{(L, \phi[d',u',x])}}}$$
is commutative, so that
\begin{eqnarray*}
\psi_d(u) &=& \|u\| \psi_{(L, \phi[d,\phi_d, u,x])}(1)\\
&=& \frac{\|u\|\|u'\|}{\|u\|} \psi_{(L, \phi[d,\phi_d, u,x])}(\frac{\|u\|}{\|u'\|})\\
&=& \frac{\|u\|\|u'\|}{\|u\|} \psi_{(L, \phi[d,\phi_d, u,x])}(\beta(1))\\
&=&\|u'\| \psi_{(L, \phi[d', \phi_{d'}, u', x])}(1)\\
&=& \psi_{d'}(u').\end{eqnarray*}
Thus,
$$T(x)= \psi_d(u) = \|u\|\psi_{(L, \phi[d,\phi_d, u,x])}(1)$$ is a judicious choice. Or would be, were this $T$ a linear contractive operator. The contractive part is clear since, no matter what it seems,
$$\|T(x)\| \leq \inf \{\|u\|: x=\varphi_d(u)\}\leq \|x\|.$$
The map $T$ is homogeneous since $T(\lambda x)= \psi_d(\lambda u) = \lambda \varphi_d(u) = \lambda T(x)$, so it only remains the additivity: why $T(x+y)=Tx + Ty$? Well, because given two independent $x, y$ (otherwise...) then we also have the element $(\ell_1^2, \varphi)$ with $\varphi: \ell_1^2 \to X$ given by $\varphi(\lambda e_1 + \mu e_2)=\lambda \frac{x}{\|x\|} + \mu \frac{y}{\|y\|}$ that, together with the canonical inclusions $\imath(1)=e_1$, $\jmath(1)=e_2$ form a commutative diagram
$$\xymatrix{
L\ar[rr]^{\imath}\ar[drr]_{ \varphi \imath}&& \ell_1^2\ar[d]^{\varphi}&& L\ar[ll]_{\jmath}\ar[dll]^{\varphi \jmath}\\
&& X &}$$
If we call $d_1= (\ell_1^1=[x], \varphi_1)$, where $\varphi_1: \ell_1^1=[x] \To X$ is given by $\varphi_1(x)=x$ and $d_2=(\ell_1^1=[y], \varphi_2)$ with $\varphi_2:\ell_1^1=[y]\To X$ given by $\varphi_2(y) = y$, then $\varphi \imath = \phi[d_1, \varphi\imath, x, x]$ and $\varphi \jmath = \phi[d_2, \varphi\jmath, y, x]$. Thus, the diagram above becomes
$$\xymatrix{
L\ar[rr]^{\imath}\ar[drr]_{ \phi[1,\varphi\imath, x,x]}&& \ell_1^2\ar[d]^{\varphi}&& L\ar[ll]_{\jmath}\ar[dll]^{\phi[2,\varphi\jmath, y,y]}\\
&& X &}$$
then, by commutativity, $\psi_{(\ell_1^2, \varphi)}\imath = \psi_{(L, \phi[1,x,x])}$
and $\psi_{(\ell_1^2, \varphi)}\jmath = \psi_{(L, \phi[2,y,y])}$. Let us assign the index $D$ to $(\ell_1^2, \varphi)$ and call $w=\|x\|e_1 + \|y\|e_2$ to write
\begin{eqnarray*}
T(x+y)&=&  \|w\| \psi_{(L, \phi[D,\varphi, w,x+y])}(1)\\
T(x)&=&  \|x\| \psi_{(L, \phi[d_1,\varphi \imath, x,x])}(1)\\
T(y)&=&  \|y\| \psi_{(L, \phi[d_2,\varphi \jmath, y,y])}(1).\end{eqnarray*}
Now, $\psi_D \imath= \psi_{d_1}$ and $\psi_D \jmath= \psi_{d_2}$ by commutativity
and thus
\begin{eqnarray*}
T(x+y)&=& \|w\| \psi_{(L, \phi[D,\varphi, w,x+y])}(1)\\
& = &\psi_D(w)\\
&=& \psi_D(\|x\| e_1 + \|y\| e_2) \\
&=& \|x\| \psi_{(L, \phi[1,\varphi \imath, x,x])} (1) + \|y\| \psi_{(L, \phi[1,\varphi \jmath, y,y])} (1)\\
& =& T(x) + T(y).\qedhere\end{eqnarray*}
\end{proof}

And, consequently,
\begin{proposition} Every dual Banach space $X^*$ is a projective limit $X^*=\lim\limits_{\leftarrow} \ell_\infty^n$ in $\mathbf{Ban_1}$.
\end{proposition}

There are different ways to prove this, being the shortest one to appeal to $\left( \lim\limits_{\rightarrow} X_i\right)^* = \lim\limits_{\leftarrow} X_i^*$, something that, at least for finite dimensional $X_i$ anyone can do by hand. A more sophisticated way could be to recall from \cite[Proposition 5.5]{HHI} that contravariant functors adjoint on the right transform inductive limits into projective limits
to then recall from \cite[Proposition 5.7]{HHI} that the duality functor is adjoint to itself on the right.

Some of the most elegant applications of comma categories to Banach space theory, the construction of universal operators, were obtained by Kubi\'s and his collaborators \cite{kub1,kub3} and have been no less elegantly exposed in HMBST, Section 6.4. See also \cite{kub2}.

\section{Kan extensions}

We are ready now to know what Kan extensions are. Extension theorems are central in mathematics, and twicely central in functional analysis. One has a mathematical object $C$ ``inside" another $D$ and an arrow or function $f: C\to E$ and wants to obtain an extension $F:D\to E$. Even when this can be done (Tietze's theorem, Hahn-Banach theorem,...) the extension is not, as a rule, unique and the way of performing the extension is not canonical. Quoting MacLane \cite{mac} ``However, if $\CAT C$ is a subcategory of $\CAT D$ each functor $\mathcal T: \CAT C\To \CAT E$ has in principle two canonical (or extreme) ``extensions" to functors $\CAT D \To \CAT E$". We will call $L\mathcal T$ and $R \mathcal T$ those extensions.\\

To see how to proceed, observe that given a covariant functor $\delta: \CAT C \To \CAT D$ between two categories and another category $\CAT E$ there is a natural functor $\delta^*: \CAT E^{\CAT D}\To \CAT E^{\CAT C}$ which is just plain composition with $\delta$; i.e., it associates to each functor $\mathcal S$ the composition functor $\mathcal S \delta$. This functor $\delta^*$ may have or have not a left adjoint $L\xdash \delta^*$ or a right adjoint $\delta^*\xdash R$. If it does,
the functors $L, R$ provide apparently good natured extensions of functors $\CAT C\To \CAT E$ to functors $\CAT D\To \CAT E$ through $\delta$.

\subsection{Left Kan extension}

\adef The Left Kan extension of $\mathcal F$ through $\delta: \CAT C \To \CAT D$ is a functor $L\mathcal F: \CAT D\To \CAT E$ together with a natural transformation $\varepsilon: \mathcal F \To L\mathcal F \delta$, with the property of being universal for this diagram; namely, given a functor $\mathcal G:\CAT D\To \CAT E$ and a natural transformation $\eta: \mathcal F \To \mathcal G \delta$ there is a unique natural transformation $\sigma: L\mathcal F\To \mathcal G$ so that $\eta = \sigma \varepsilon $.\zdef

To show its distinguished character, the functor $L\mathcal F: \CAT D \To \CAT E$ comes accompanied with a natural transformation providing its uniqueness property. Let us draw the situation: this is the extension functor
$$\xymatrix{\CAT C\ar[r]^\delta \ar[rd]_{\mathcal F}& \CAT D \ar[d]^{L\mathcal F}\\
&\CAT E}$$
and here it is the associated natural transformation
\begin{equation}\label{LKan}\xymatrix{\CAT C \ar[dd]_{\mathcal F}\ar@{=}[rr]& &\CAT C \ar[rr]^\delta \ar[dd]^{L\mathcal F\delta}&&\CAT D \ar[dd]^{L\mathcal F}\\
\ar[rr]_{\varepsilon}&&\\
\CAT E \ar@{=}[rr] &&\CAT E \ar@{=}[rr]& &\CAT E}\end{equation}
and here is the uniqueness
$$\xymatrix{L\mathcal F \ar[rr]_{\sigma}& &\mathcal G\\
L\mathcal F\delta \ar[rr]_{\sigma\delta}& &\mathcal G\delta \\
&\mathcal F \ar[ul]^{\varepsilon}\ar[ur]_{\eta}&}$$
One has:

\begin{proposition}\label{leftKan} If $\delta^*$ has a left adjoint $L\xdash \delta^* $ then
$L\mathcal F$ is the Left Kan extension of $\mathcal F$. Thus, the Left Kan extension of $\mathcal F$ comes defined by the equation $[L\mathcal F, \mathcal G] = [\mathcal F, \mathcal G\delta]$ for every $\mathcal G:\CAT D\To \CAT E$.\end{proposition}

This proof we will make in the next section for the Right Kan extension; so the reader just have to dualize it. If $\delta^*$ has a left adjoint then every functor has a Left Kan extension, but it could happen that a given functor admits Left Kan extension even when $\delta^*$ has no left adjoint. The next method shows a way to compute $L\mathcal F$.

\begin{proposition}\label{leftKan2} If the category $\CAT E$ admits limits then
 $$L\mathcal F(D)= \lim \limits_{\longrightarrow} \mathcal F^D$$
\end{proposition}
\begin{proof} Let $ \mathcal F: \CAT C\To \CAT E$ be a functor. Pick an object $D$ of $\CAT D$, form the comma category $(\CAT C \Rsh D)$ and then the induced functor $\mathcal F^D: (\CAT C \Rsh D) \To \CAT E$. The limit $\lim \limits_{\longrightarrow} \mathcal F^D$
exists by hypothesis and we can form the functor $L\mathcal F$ as above. Our plan is to show that it defines the Left Kan extension of $\mathcal F$. It is clear that $L\mathcal F(C) = \lim \limits_{\longrightarrow} \mathcal F^C = \mathcal  FC$ for every object $C$ of $\CAT C$. Thus $L\mathcal F \delta =\mathcal F $ and we can set $\varepsilon = 1$ as the accompanying natural transformation $\varepsilon: \mathcal F \To L \mathcal F \delta$. Let $\mathcal G$ be another functor $\CAT D\To \CAT E$ with accompanying natural transformation $\eta: \mathcal F\To \mathcal G\delta$. Namely, for given $c: C\To C'$ there is a commutative diagram
$$\xymatrix{
\mathcal F(C)\ar[d]_{\mathcal F(c)}\ar[r]^-{\eta_C}& \mathcal G(C)\ar[d]^{\mathcal G(c)}\\
\mathcal F(C')\ar[r]_-{\eta_{C'}}& \mathcal G(C')
}$$
It is then guaranteed the existence of the natural transformation $\sigma: L\mathcal F\To \mathcal G$ we need, namely that for given $d: D\To D'$ one has a commutative diagram
$$\xymatrix{
\lim\limits_{\longrightarrow} \mathcal F^D \ar[d]_{L\mathcal F(d)}\ar[r]^-{\sigma_D}& \mathcal G(D)\ar[d]^{\mathcal G(d)}\\
\lim\limits_{\longrightarrow} \mathcal F^{D'} \ar[r]_-{\sigma_{D'}}& \mathcal G(D')
}$$
To make appear the arrow $\sigma_D$ recall that objects of the comma category $(\CAT C \Rsh D)$ are pairs $(C,f)$
with $C$ an object of $\CAT C$ and $f: C\to D$ an arrow. There exists therefore arrows $G(f): G(C)\To G(D)$ and thus the universal property of limits yields $\sigma_D$.\end{proof}

\subsection{Right Kan extension}

\adef The Right Kan extension of $\mathcal F$ through $\delta: \CAT C \To \CAT D$ is a functor $R\mathcal F: \CAT D\To \CAT E$ together with a natural transformation $\varepsilon: R\mathcal F \delta\To \mathcal F$, with the property of being universal for this diagram; namely, given a functor $\mathcal G:\CAT D\To \CAT E$ and a natural transformation $\eta: \mathcal G \delta\To \mathcal F$ there is a unique natural transformation $\sigma: \mathcal G\To R\mathcal F$ so that $\eta = \varepsilon \sigma\delta$.\zdef

The existence of natural transformation accompanying the extension that enjoys some uniqueness property
provides a clean way of saying that the functor $R\mathcal F: \CAT D \To \CAT E$ is a distinguished extension of $\mathcal F$. Let us draw the required properties: here it is the extension functor
$$\xymatrix{\CAT C\ar[r]^\delta \ar[rd]_{\mathcal F}& \CAT D \ar[d]^{R\mathcal F}\\
&\CAT E}$$
here it is the associated natural transformation
\begin{equation}\label{RKan}\xymatrix{\CAT C \ar[dd]_{\mathcal F}\ar@{=}[rr]& &\CAT C \ar[rr]^\delta \ar[dd]^{R\mathcal F\delta}&&\CAT D \ar[dd]^{R\mathcal F}\\
&&\ar[ll]_{\varepsilon}\\
\CAT E \ar@{=}[rr] &&\CAT E \ar@{=}[rr]& &\CAT E}\end{equation}
and here it is the uniqueness
$$\xymatrix{R\mathcal F & &\mathcal G\ar[ll]_{\sigma}\\
R\mathcal F\delta \ar[dr]_{\varepsilon}& &\mathcal G\delta \ar[dl]^{\eta}\ar[ll]_{\sigma\delta}\\
&\mathcal F &}$$

One has the dual version of Proposition \ref{leftKan}, now with proof:
\begin{proposition}\label{rightKan} If $\delta^*$ has a right adjoint $\delta^* \xdash R$ then
$R\mathcal F$ is the Right Kan extension of $\mathcal F$. Thus, the Right Kan extension of $\mathcal F$ comes defined by the equation $[\mathcal  G, R\mathcal  F] = [\mathcal G\delta, \mathcal F]$ for every $\mathcal G:\CAT D\To \CAT E$.\end{proposition}
\begin{proof} Set $\mathcal G=R\mathcal F$ in the equation $[\mathcal G, R\mathcal F] = [\mathcal G\delta, \mathcal F]$ that thus becomes  $[R\mathcal F, R\mathcal F] = [R\mathcal F\delta, \mathcal F]$ and pick the natural transformation $\varepsilon: R\mathcal F\delta\to \mathcal F$ corresponding to the identity $R\mathcal F\to R\mathcal F$.
Now, if $\mathcal G$ is another functor accompanied with a natural transformation $\eta\in [\mathcal G\delta, \mathcal F]$ then obtain the corresponding $\sigma \in [\mathcal G, R\mathcal F]$. To check that $\varepsilon \sigma\delta=\eta$ as elements of $[\mathcal G\delta, \mathcal F]$ one just has to check that the corresponding elements in $[\mathcal G,R\mathcal F]$ are the same: $\sigma$. And this is simple since $\sigma\delta\in[\mathcal G\delta, R\mathcal F\delta]=[\mathcal G, R(R\mathcal F\delta)]=[\mathcal G, R\mathcal F]$ corresponds to $\sigma$
then $\varepsilon \sigma\delta$ corresponds to $1\sigma= \sigma$.\end{proof}

Of course that if $\delta^*$ has a right adjoint then every functor has a Right Kan extension, but it could happen that a given functor admits Right Kan extension even when $\delta^*$ has no right adjoint. The following characterization and its proof are clear dualizations of
Proposition \ref{leftKan2}

\begin{proposition}\label{raghtKan2} If the category $\CAT E$ admits colimits then
 $$R\mathcal F(D)= \lim \limits_{\longleftarrow} \mathcal F_D$$
\end{proposition}

\section{Uses of Kan extensions for we, Banach spacers}\label{forwe}

The category of Banach spaces does not live alone in Banach space affairs: it is intimately connected with the categories $\mathbf{pBan}$ of $p$-Banach spaces and with the more obscure one $\mathbf {QBan}$ of quasi Banach spaces. Now, $\mathbf {Ban}$ is not, categorically speaking, a good category because operators do not have the right cokernels, arbitrary products do not exist, let alone limits... However, it has, at least, enough injective and projective objects (in a slightly more restrictive sense \cite{potho}: after all, the operators one can actually deal with are those with closed range (see Section \ref{heart}), and where \emph{enough} means enough to \emph{present} any Banach space). The category $\mathbf {pBan}$ is worse: in addition to that, it has no injectives at all \cite{hmbst}. Luckily, it still has (in the same restricted sense) enough projectives. But, alas! There is no minimum in the slope of worst, and $\mathbf{ QBan}$ is awful: in addition to everything that was already wrong it does not have injective or projective objects. How operative $\mathbf{Ban}$ is is a returning topic in this paper. Let us focus now in how bad $\mathbf{QBan}$ is.

Moving at bullet-time, skipping next section and arriving to Section \ref{sec:derivation} we will encounter that there are two standard procedures to proceed with the derivation of functors: one uses injective objects and yields right-derivation and the other uses projective  objects and yields left-derivation. The good news is that none of them is available in $\mathbf {QBan}$, so we need to rethink derivation in this category from scratch. In $\mathbf {pBan}$ one at least can still perform left-derivation and that fosters the delusion of understanding. But not in $\mathbf {QBan}$. We could think this is important because Homology is based on the idea of derivation (in approximately the same sense that Calculus is based on the idea of derivative). But we could also adopt a broader view, not stick to derivation, and consider that Homology is based on the idea that one can obtain an algebraic map of topological spaces. In which case to put in focus $\mathbf {QBan}$ seems yet more important. To provide us with a workable slogan: Is $\Ext$ the derived functor of $\mathfrak L$ in $\mathbf {QBan}$?\\

Well, the first issue is: which functors are exactly $\mathfrak L$ and $\Ext$? Because the definition of a functor requires to specify the categories between which it acts. And $\mathfrak L$ is deceitful. And it is deceitful precisely because it is a functor $\mathbf{ Ban}\To \mathbf{Ban}$; thus, if we trust that $\Ext$ is going to be the derived functor of $\mathfrak L$ then we must find a way to make
$\Ext$ a functor $\mathbf{ Ban}\To \mathbf{Ban}$. This, to the best of our current knowledge, is impossible: $\Ext(\ell_2, \ell_2)$ is not and cannot be a Banach space, see HMBST, Section 4.5: it is not because the ``natural norm" to be set on $\Ext(\ell_2, \ell_2)$, which is moreover the one that makes continuous the maps appearing in the homology sequence, is not Hausdorff; namely, there are $0\neq \alpha  \in \Ext(\ell_2, \ell_2)$ for which $\|\alpha\|=0$ (HMBST, Proposition 4.5.5).  Mathematical life is funny here, because when $p<1$ the space $\Ext(L_p/\ell_2, \R)$ is a Hilbert space, (HMBST, Corollary 4.5.2). The part of me that believes that Tyrion was a rightful candidate to rule the seven kingdoms grunts that having or not having a trivial dual cannot be the responsible for the difference of behaviours because $\Ext(\ell_2, \ell_2)$ is the same in $\CAT{Ban}$ than in $\CAT{QBan}$. While the part of me that believes that it should have been Gandalf who sits in the iron throne
is convinced that it is. Be as it may, and summing up: we cannot support or endorse in any natural way that $\Ext$ is the derived functor of the functor $\mathfrak L: \mathbf{ Ban}\To \mathbf{Ban}$.

What we can do and do is to slip away and consider the functor $\mathfrak L: \mathbf{ Ban}\To \mathbf{Vect}$ (the category of vector spaces and linear maps). Now we can prove that $\Ext$ is the derived functor of this $\mathfrak L$! And who could be the derived functor or
$\mathfrak L: \mathbf{ QBan}\To \mathbf{Vect}$? Well, $\Ext$ behaves as if it were: its construction can be done Yoneda style, as well as the iterated $\Ext^n$ functors and the long homology sequence that connects them exists in complete analogy with the $\mathbf{Ban}$ case. Is this enough to accept that $\Ext$ is \emph{the} derived functor of $\mathfrak L$ in $\mathbf{QBan}$ through a derivation process that could not exist? Of course not. But the pot of gold at the end of this rainbow is whether, and how, derivation can be defined in $\mathbf{QBan}$.

One idea we can explore in this section is to what extent $\Ext_{\mathbf{QBan}}$ is \emph{a} Kan extension of $\Ext_{\mathbf{Ban}}$. This would at least provide us some basis to believe that a Kan extension of a derived functor can be called a derived functor. Moving in that direction, let us gain first some confidence in obtaining Kan extensions of Banach functors. Observe that since the Semadeni-Zidenberg theorem \cite{semazide} establishes that $\mathbf {Ban_1}$ admits limits and colimits, one has:

\begin{cor} Let $\mathcal F:\mathbf {Ban_1}\To \mathbf{Ban_1}$ be a Banach functor.
\begin{itemize}
\item If $R\mathcal F$ is a Right Kan extension
through $\square$ to $\mathbf{Q Ban_1}$ then $R\mathcal F\square=\mathcal F$
\item If $L\mathcal F$ is a Left Kan extension through $\square$ to $\mathbf{QBan_1}$ then $L\mathcal F\square=\mathcal F$.
\end{itemize}
\end{cor}

Let us now test our abilities obtaining the Kan extensions $\mathbf {QBan_1}\To \mathbf {Ban_1}$ of the identity $1: \mathbf {Ban_1}\To \mathbf{Ban_1}$:

\begin{proposition}\label{R1L1} $L1 = R1 = co$.\end{proposition}
\begin{proof} It is obvious that $co \square =1$ and thus the natural transformation $\varepsilon: co \square \To 1$ we st is the identity. To show that $co$ is the Right Kan extension of $1$, let us assume first that one has another functor $\mathcal G:\mathbf {QBan_1}\To
\mathbf {Ban_1}$ so that $\mathcal G\square = 1$, for which the natural transformation $\eta $ is also the identity. To obtain the required natural transformation $\sigma: \mathcal G\To co$, pick a quasi Banach space $Q$, consider the canonical operator $\delta_Q: Q\to co(Q)$
and form $\mathcal G(\delta_Q): \mathcal G(Q)\to \mathcal G(co(Q)) =co (Q)$. Set $\sigma_Q= \mathcal G(\delta_Q)$. In the general case, if we just assume that there is a natural transformation $\eta: \mathcal G\square\to 1$ then set $\sigma_Q = \eta_{co(Q)}\mathcal G(\delta_Q)$. This shows that $co$ is the right Kan extension of the identity. To show it is also the left Kan extension, fix a quasi Banach space $Q$ and let us show that $\lim\limits_{\rightarrow} 1^Q$ of the translation $1^Q$ of the identity functor to the comma category  $(\mathbf {Ban_1} \Rsh Q)$ is $co(Q)$: consider the family of Banach spaces $X$ and contractive operators $\tau: X\To Q$; so $B_X\subset \tau^{-1}(B_Q)$. This family includes the following elements: pick $x\in B_Q$
and set the pair $([x], \imath_x)$ where $\imath_x(x)=x$. The existence of a contractive operator $\ell_x: ([x], \imath_x)\To L$ implies $\|\ell_x(x)\|\leq \|x\|\leq 1$. This yields the existence of a contractive operator $Q \To \lim 1^Q$. Now, given a Banach space $Y$ and a contractive operator $\tau: Q\To Y$, the composition $\tau\varphi_X$ with $\varphi_X: X\To Q$ provides contractive operators $X\To Y$, and consequently a contractive operator $\lim\limits_{\rightarrow} 1^Q\To Y$. All this means that $\lim\limits_{\rightarrow} 1^Q = co(Q)$.\end{proof}

The following two diagrams are to illustrate the similitude between the two Kan extensions (see also \cite[Theorem 2.6]{lehner}). Compare to diagrams (\ref{RKan}) and (\ref{LKan}) in the previous Section:

$$\xymatrixrowsep{0,5cm}\xymatrixcolsep{1,5cm}
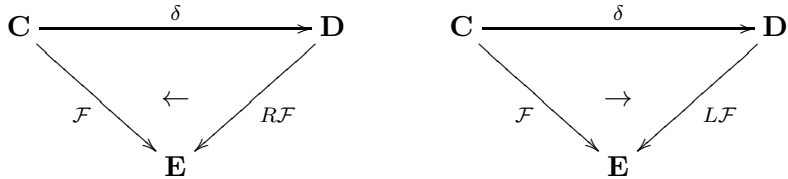
\xymatrix{
\CAT C\ar[rr]^\delta\ar[ddr]_{\mathcal F}  & & \CAT D\ar[ddl]^{R\mathcal F}\\
&\leftarrow&\\
&\CAT E&}\quad\quad\quad\xymatrix{
\CAT C\ar[rr]^\delta\ar[ddr]_{\mathcal F}  & & \CAT D\ar[ddl]^{L\mathcal F}\\
&\rightarrow&\\
&\CAT E&}
$$

Keep also in mind that a Kan extension is nowhere claimed to be \emph{an extension}, in the sense that $R\mathcal F\delta =\mathcal F$ (when this makes sense, as it occurs with our $\square: \mathbf {Ban} \To \mathbf {QBan}$ case. Namely, even if $\mathcal G: \mathbf {QBan} \To \mathbf {Ban}$ is a functor such that $\mathcal G\square = \mathcal F$ that does not guarantee that $\mathcal G$ is the Right Kan extension of $\mathcal F$ because there is no way to guarantee the existence of the natural transformation $\sigma: \mathcal G \To R\mathcal F$. \\

Time to splash. What we did in Proposition \ref{R1L1} was by no means showing that $\square^*$ has adjoints, but to calculate that $R1=L1 = co$, that is, to showing that one has the identities $[\mathcal G\square, 1] = [\mathcal G, co]$ and $[co, \mathcal G] = [1, \mathcal G\square]$ for every $\mathcal G: \mathbf {QBan}\To \mathbf {Ban}$. And those two things, getting an adjoint for $\square^*$ and getting some identities, are
different things, by far. Well, not that far. Let us prepare the way:

\begin{lemma}\label{vv*} Given functors $\mathcal U: \CAT C \To \CAT D$ and $\mathcal V: \CAT D \To \CAT C$, if $\mathcal V \xdash \mathcal U$ then $\mathcal U^* \xdash \mathcal V^*.$
\end{lemma}
\begin{proof} One claims that given $\mathcal G\in \CAT E^{\CAT D}$ and $\mathcal F\in \CAT E^{\CAT C}$ then $[\mathcal U^*\mathcal G, \mathcal F] = [\mathcal G, \mathcal V^* \mathcal F]$ which is as obvious as obscure: given morphisms $c: C\To C'$ and $d: D\To D'$ one claims that
natural transformations $\eta, \nu$ as in
$$\xymatrix{
\mathcal G\mathcal U C\ar[d]_{\mathcal G \mathcal Uc}\ar[r]^{\eta_C}& \mathcal FC\ar[d]^{\mathcal Fc}\\
\mathcal G \mathcal U C'\ar[r]_{\eta_{C'}}& \mathcal F C' }\quad \quad \quad \quad \xymatrix{
\mathcal G D\ar[d]_{\mathcal G d}\ar[r]^{\nu_D}& \mathcal F \mathcal V D\ar[d]^{\mathcal F \mathcal V g}\\
\mathcal G D' \ar[r]_{\nu_{D'}}& \mathcal F  \mathcal V D'}$$
correspond one to another, as it is clear they do (recall there are natural transformations $\varepsilon: 1\to \mathcal U \mathcal V$ and $\varepsilon': \mathcal V\mathcal U\to 1$): set  $\nu_D = \eta_{\mathcal V D}$ to  get
$$\xymatrix{
\mathcal G\mathcal U \mathcal V D\ar[d]_{\mathcal G \mathcal U \mathcal V d}\ar[r]^{\eta_C}& \mathcal F \mathcal V D\ar[d]^{\mathcal F \mathcal V d}\\
\mathcal G \mathcal U \mathcal V D'\ar[r]_{\eta_{\mathcal V D'}}& \mathcal F \mathcal V D' }
\Longrightarrow \xymatrix{
\mathcal G D \ar[d]_{\mathcal Gd}\ar[r]^-{\varepsilon_D}&\mathcal G \mathcal U \mathcal V D\ar[d]_{\mathcal G \mathcal U \mathcal V d}\ar[r]^{\eta_{\mathcal V D}}& \mathcal F \mathcal V D\ar[d]^{\mathcal F \mathcal V g}\\
\mathcal G D' \ar[r]_-{\varepsilon_{D'}}& \mathcal G \mathcal U \mathcal V D' \ar[r]_{\eta_{\mathcal V D'}}& \mathcal F  \mathcal V D'}$$
and analogously in the other case: set $\eta_{C} = \nu_{\mathcal U C}$ and use $\varepsilon'$.\end{proof}

In particular (and probably simpler to prove) we get $\square^*\xdash co^*$ as it follows from $co\;\xdash \square$. Thus, by the uniqueness of adjoints, one gets:

\begin{theorem}\label{righthom} $co^* = R$. Thus, given a Banach functor $\mathcal F$, or a functor $\mathcal F: \mathbf{Ban}\To \mathbf{QBan}$ or $\mathcal F: \mathbf{Ban}\To \mathbf{Vect}$, its Right Kan extension is $co^* \mathcal F$.\end{theorem}

In particular: let $X$ be a Banach space and let $Q', Q$ be quasi Banach spaces:
\begin{itemize}
\item For $\mathfrak L_{Q'}: \mathbf{Ban}\To \mathbf{Ban}$ or  $\mathfrak L_{Q'}: \mathbf{Ban}\To \mathbf{QBan}$
one has $$R\mathfrak L_{Q'} (Q) = \mathfrak L(Q', co(Q)).$$
\item For $\mathfrak L^{X}: \mathbf{Ban}\To \mathbf{Ban}$ or  $\mathfrak L^{Q'}: \mathbf{Ban}\To \mathbf{QBan}$
one has $$R\mathfrak L^X (Q) = \mathfrak L(co(Q), X) = \mathfrak L(Q, X) \quad \quad \mathrm{and} \quad\quad R\mathfrak L^{Q'}(Q) = \mathfrak L (co(Q), Q').$$
\item For $\Ext_{\mathbf{Ban}}(X, \cdot): \mathbf{Ban}\To \mathbf{Vect}$ one has $R\Ext_{\mathbf{Ban}}(X, \cdot)(Q) = \Ext_{\mathbf{Ban}}(X, co(Q))$.
\item For $\Ext_{\mathbf{QBan}}(Q', \cdot): \mathbf{Ban}\To \mathbf{Vect}$ one has $R\Ext_{\mathbf{QBan}}(Q', \cdot)(Q) = \Ext_{\mathbf{QBan}}(Q', co(Q))$.
\item For $\Ext_{\mathbf{Ban}}(\cdot, X): \mathbf{Ban}\To \mathbf{Vect}$ one has $R\Ext_{\mathbf{Ban}}(\cdot, X)(Q) = \Ext_{\mathbf{Ban}}(co(Q), X)$.
\item For $\Ext_{\mathbf{QBan}}(\cdot, Q'): \mathbf{Ban}\To \mathbf{Vect}$ one has $R\Ext_{\mathbf{Ban}}(\cdot, Q')(Q) = \Ext_{\mathbf{Ban}}(co(Q), Q')$.
\end{itemize}
What about the left Kan extensions of $\mathfrak L_X, \mathfrak L^X$? Well, a very much expected surprise is that our first hazardous guess $L=R$ cannot actually be: $L\xdash \square^*$ and $L=R$ yields $co^*\xdash \square$, and this is at a teeth skin from saying $\square\;\xdash co$, which is false. And false on its own, anyway. By the pointwise formula in Proposition \ref{leftKan2} we get:
\begin{theorem}\label{lefthom} Let $X$ be a Banach space and $Q$ a quasi Banach space. Given $\mathfrak L_X$ taking values in
either $\mathbf{Vect}$ or $\mathbf{Ban_1}$ one has
\begin{itemize}
\item $L \mathfrak L_X(Q)= \lim \limits_{\longrightarrow} (\mathfrak L_X)^Q$
\item $L\mathfrak L^X (Q) = \lim \limits_{\longrightarrow} (\mathfrak L^X)^Q$.
\end{itemize}
Given $\Ext_{\mathbf B}(X, \cdot)$ taking values in $\mathbf {Vect}$ one has
\begin{itemize}
\item $L\Ext_{\mathbf{Ban}}(X, \cdot)(Q) = \lim \limits_{\longrightarrow} \Ext_{\mathbf{Ban}}(X, \cdot)^Q$.
\item $L\Ext_{\mathbf{Ban}}(\cdot, X)(Q) = \lim \limits_{\longrightarrow} \Ext_{\mathbf{Ban}}(\cdot, X)^Q$
\end{itemize}
\end{theorem}
Always taking advantage from the fact that the categories $\mathbf{Vect}$ and $\mathbf{Ban_1}$ have limits and colimits. We have not a plausible guess about what occurs with those functors when taking values in $\mathbf {QBan_1}$ since the Semadeni-Zidenberg theorem does not work in $\mathbf {QBan_1}$ and the limit formulae are no longer available.\\

Is this what we expected? Likely not: one would have wanted, perhaps, that the true natural extension of $\mathfrak L$ to be the same $\mathfrak L$. But this is not how things are, as it is demonstrated by the work of Herz and Pelletier \cite{hpelle} and then Pelletier \cite{pelle1,pelle2}. Inspired, to some extent, by the work of Cigler \cite{cig} as it is mentioned in \cite[p.486]{pelle1}, they show that if $\CAT B$ is a suitable subcategory of $\mathbf {Ban_1}$ then the Left Kan extension of the functor $\mathfrak L$ to $\mathbf {Ban_1}$ is formed by the operators that factorize through an element of $\CAT B$. In particular, if $\CAT R$ is the category of reflexive spaces, its Left Kan extension are the weakly compact operators. And if $\CAT F$ is the subcategory of finite dimensional spaces its Left Kan extension is, under some Approximation Property hypotheses, the compact operators. Very nice. And showing that \emph{the} extension of $\mathfrak L$ should not be expected to be $\mathfrak L$. Let alone that the natural extension of $\Ext_{\mathbf{Ban}}$ should be $\Ext_{\mathbf{QBan}}$. Oh, well, worse things happen at sea.

\section{Homology}\label{homo}
After 540 pages of HMBST we (authors, readers, by-passers skimming the pages....?) should have formed an idea about what Banach space homology is and what is it for. The fulcrum point is the idea hinted above that a mathematical
object is the object \emph{and} its many \emph{representations}, whatever that might mean. An analogue could be the (quite standard) idea that a real number $r$ is not only the number $r$ but also all sequences of rational numbers converging to $r$. In categorical jargon, the category of real numbers is derived from the category of rational numbers by: fixing as objects (certain) sequences of rational numbers plus an equivalence relation that suitably transforms morphisms of sequences into morphisms in the derived category. To do this uplifting from, say, $\mathbf{Ban}$ to $\mathbf{Ban}^{\mathbf{Ban}}$ we need the core ideas of ``complex" and ``exact sequence":

\adef\label{def:sex} An exact sequence of quasi Banach spaces is a diagram formed by quasi Banach
spaces and operators $\xymatrix{
\cdots \ar[r]&X_{i-1}\ar[r] & X_i  \ar[r]& X_{i+1}\ar[r] &\cdots  }$ in which the kernel of each arrow coincides with the image of the
preceding one.\zdef

A morphism between two exact sequences is obviously a sequence $(f_i)$ of morphisms $f_i: X_i\to Y_i$ yielding commutative diagrams
$$\xymatrix{
X_i\ar[r]\ar[d]_{f_i}& X_{i+1}\ar[d]^{f_{i+1}}\\
Y_i\ar[r]& Y_{i+1}}$$

Let us call $\mathbb{EX}$ the category of exact sequences of Banach spaces. An especially distinguished type of exact sequences is that of \emph{short} exact sequences
$$
\xymatrix{
0 \ar[r]&Y \ar[r] & Z  \ar[r]& X\ar[r] & 0  }
$$

The arrows in the category are provided by triples of arrows $(\alpha, \beta, \gamma)$ making commutative diagrams
\begin{equation*}\label{triple1}\xymatrix{0\ar[r]& \cdot \ar[r]\ar[d]_\alpha& \cdot \ar[r]\ar[d]^\beta &  \cdot \ar[r]\ar[d]^\gamma & 0\\
0\ar[r]& \cdot \ar[r]& \cdot \ar[r] &  \cdot \ar[r] & 0}\end{equation*}

And we thus arrive to the category $\mathbb{S}$ of short exact sequences. Two short sequences are thus isomorphic when there is a diagram as above with $\alpha, \beta, \gamma$ isomorphisms, which is the notion of \emph{isomorphic} sequences used in HMBST. Good examples of (long) exact sequences in $\mathbf {Ban}$ are provided by the so-called \emph{projective presentations} of a given space $X$; namely, exact sequences
$$\begin{CD} \cdots @>>> P^{-n} @>>> \cdots @>d^2>> P^{-2} @>{d^{1}}>> P^{-1} @>{d^0}>> X @>>>0\end{CD} $$
in which each space $P^{-n}$ is projective (i.e., some space $\ell_1(\Gamma)$, since these are the only projective (in our restricted sense) Banach spaces \cite{kothe} (see also \cite{groth} for a close result); or else \cite{orto} where it is proved that $\ell_p(I)$ are the only projective (in restricted sense) $p$-Banach spaces \cite{orto}. Alternatively, \emph{injective presentations} namely, exact sequences
$$\begin{CD} 0@>>>X @>{\partial^0}>> I^{1} @>{\partial^1}>> I^2 @>>>\cdots @>>> I^n @>>>\cdots\end{CD} $$
in which each space $I^{n}$ is injective.  Examples of short exact sequences are provided by a Banach space $Z$, a subspace $Y$ and the corresponding quotient $Z/Y$ in the form
$$
\xymatrix{
0 \ar[r]&Y \ar[r] & Z  \ar[r]& Z/Y\ar[r] & 0  }
$$
Every short exact sequence of Banach spaces is isomorphic to one of this type. But beyond The Exactness Wall there is much more life. The wildlings in this case are \emph{the complexes}. A complex is, well, the ``non exact" part of an exact sequence; that is, a sequence
$\xymatrix{
\cdots \ar[r]&X_{i-1}\ar[r] & X_i  \ar[r]& X_{i+1}\ar[r] &\cdots  }$ such that the composition of any two consecutive arrows is $0$.

So the reader may ask: why do we bother with complexes when we can have exact sequences? There are various reasons. One is that the notion of complex makes sense in every category having $0$ object while that of exact sequence is much more delicate.\\

\textbf{Intermission about the future.} In a category such as Abelian groups, an exact sequence is a complex in which the kernel of each arrow coincides with the image of the preceding. The notion of kernel passes without difficulties to an arbitrary category but that of ``image" does not. One has the cokernel notion, of course, and probably one would like to say that a sequence is exact (at least short exact) when each arrow is the kernel of the next one and the cokernel of the previous one. In an Abelian category (see later) things go more or less straight since one can define (see Definition \ref{image}) the image of an arrow as the kernel of its cokernel (see \cite[p.6]{weibel}) and proceed. But that is because an Abelian category is, by definition, a place where monics are kernels of their cokernels and epics are cokernels of their kernels. \textbf{End of the future intermission.}\\

Thus, having exact sequences requires a category where kernels and cokernels exist and behave well. It is curious that Banach spaces, which is not actually one of those places, has a clean notion of exact sequence with it (this topic shall be discussed in Section \ref{heart}).
Another reason is that, even if one has exact sequences, more general complexes could be needed (something similar to: to deal with exact sequences of Banach spaces we need quasi Banach spaces) to do things that cannot be done with only exact sequences. Moving on, the typical complex one will consider will have the form
$$\begin{CD} \cdots @>>> C^{-n} @>>> \cdots @>d^1>> C^{-1} @>{d^0}>> X @>{\partial^0}>> C^1 @>{\partial^1}>> \cdots @>>> C^n@>>>\cdots \end{CD} $$

Projective presentations are complexes as above with $C^n=0$ for all $n>0$ (and only nontrivial $d^n$ maps) and injective presentations are complexes with $C^n=0$ for all $n<0$ (and only nontrivial $\partial^n$ maps). Let us call $\mathbb{COM}(\CAT C)$ the category of complexes
of a category $\CAT C$. In particular $\mathbb{COM}(\mathbf B)$ is the category of Banach space complexes. Isomorphism in this category are what they are. Let us reconsider if that is that what we want: While working in the \emph{short} context, a projective presentation of a Banach space $X$ is a short exact sequence $0\To \kappa \To \ell_1(\Gamma) \To X \To 0$. A result that is half-classical and half-half-homological \cite[Section 2.7 and 2.11.6]{hmbst} is that if $0\To \kappa \To \ell_1(\Gamma) \To X \To 0$ and $0\To \kappa' \To \ell_1(\Gamma') \To X \To 0$
are two projective presentations of $X$ then there is a commutative diagram
\begin{equation}\label{triple1}\xymatrix{0\ar[r]& \kappa\times \ell_1(\Gamma') \ar[r]\ar[d]_\alpha& \ell_1(\Gamma)\times \ell_1(\Gamma') \ar[r]\ar[d]^\beta &  X \ar[r]\ar@{=}[d] & 0\\
0\ar[r]& \kappa'\times \ell_1(\Gamma) \ar[r]& \ell_1(\Gamma')\times \ell_1(\Gamma) \ar[r] &  X \ar[r] & 0}\end{equation}
with $\alpha, \beta$ isomorphisms. In other words, the two projective presentations are isomorphic ``up to some projective space". Thus, if we want that ``all projective presentations are the same", that is what we must get. Even if the isomorphism notion above for exact sequences is useful and has provided many deep results in Banach space theory, the right notion to work with projective presentations seems to be ``equivalent up to some projective space". Let us give shape to this.

Grasping the right notion just requires to realize that the right objects in a category of short exact sequences are \emph{equivalence} classes of exact sequences. Two short exact sequences are said to be equivalent if there exists an operator $\beta$ making a commutative diagram
\begin{equation*}\label{triple1}\xymatrix{0\ar[r]& \cdot \ar[r]\ar@{=}[d]& \cdot \ar[r]\ar[d]^\beta &  \cdot \ar[r]\ar@{=}[d] & 0\\
0\ar[r]& \cdot \ar[r]& \cdot \ar[r] &  \cdot \ar[r] & 0}\end{equation*}
So beware that equivalence is a notion defined between exact sequences having the same initial and final spaces. The equivalence class of a sequence $\texttt z$ will be denoted $[\texttt z]$. Now, triples of arrows are no longer well suited to mingle with equivalence classes of exact sequences: after all, \emph{where} should one define $\beta$? So, let us declare $\mathbb {EEX}$ to be the category having as objects equivalence classes $[\texttt z]$ of exact sequences and as morphisms couples (actually, equivalence classes of, see below) $(\alpha, \gamma)$ such that $[\alpha\,\texttt{z}]=[\texttt{z}'\gamma]$. Recall that this just means that the pushout sequence $\alpha\,\texttt{z}$ and the pullback sequence $\texttt{z}'\gamma$ are equivalent. But, acting this way, the isomorphism notion has changed to a new \texttt{isomorphism} notion that we proceed to describe. Given an exact sequence $\texttt z$
and a quasi Banach space $E$ then will denote $E\times \texttt z$ the sequence \begin{equation*}
\xymatrix{
 0\ar[r] & E\times Y  \ar[r]^-{{\bf 1}_E\times\jmath} & E\times Z \ar[r]^-{0\oplus \rho}  & X  \ar[r] & 0
}
\end{equation*}
and $\texttt z \times E$ the sequence
\begin{equation*}
\xymatrix{
 0\ar[r] &  Y  \ar[r]^-{(\jmath,\,0)} & Z\times E \ar[r]^-{\rho\times {\bf 1}_E}  & X\times E  \ar[r] & 0
}\end{equation*}
In general, these sequences will be called \emph{multiples} of $\texttt z$. Since two objects $[\texttt z]$ and $[\texttt z']$ are \texttt{isomorphic} if there exist arrows $(\alpha, \gamma): [\texttt z] \To [\texttt z'] $ and
$(\alpha', \gamma'): [\texttt z'] \To [\texttt z] $ such that $(\alpha, \gamma)(\alpha', \gamma') = 1 = (\alpha', \gamma')(\alpha, \gamma)$, it is rather clear that $[E \times \texttt z]$ and $[\texttt z \times E]$ are \texttt{isomorphic} objects (although, in
general, the sequences $E \times \texttt z$ and $\texttt z \times E$ need not be either equivalent or isomorphic (for obvious reasons!). One has \cite[Proposition 3.1]{castmoreLN}:

\begin{proposition} $ [\texttt z]$ and $ [\texttt z']$ are \texttt{isomorphic} if and only if there exist $E, E'$ so that  $E\times \texttt z \times E'$ and $E\times \texttt z' \times E'$ are isomorphic.\end{proposition}

We still need to move on to get home: to define our category $\mathbb{SEX}$ with objects as in $\mathbb{EEX}$ but with some identification of arrows: we declare $(\alpha, \gamma)$ to be $0$ when $\alpha \texttt z=0=\texttt z'\gamma$ and, consequently, $(\alpha, \gamma) = (\alpha', \gamma')$ when $[(\alpha -\alpha')\texttt z] = 0 = [\texttt z' (\gamma - \gamma')]$.
To translate this to complexes, where no natural notion of equivalence disturbs us, we will return to an equivalence-free approach. Peculiarities of short exact sequences imply that  $(\alpha, \gamma)= (\alpha', \gamma')$ in $\mathbb {SEX}$ (namely, $(\alpha - \alpha')\texttt z = 0 = \texttt z' (\gamma - \gamma')$) if and only if  there are arrows $u, v$
 \begin{equation*}
\xymatrix{ 0\ar[r] & Y \ar[r]^\jmath \ar[d]_{\alpha - \alpha'} & Z  \ar[dl]_u\ar[r]^\rho\ar[d] & X\ar[r]\ar[d]^{\gamma-\gamma'}\ar[dl]_v & 0 &(\texttt z)\\
0\ar[r] & Y' \ar[r]^{\jmath'} & Z'  \ar[r]^{\rho'} & X'\ar[r] & 0 &(\texttt z')
}
\end{equation*}
 such that $u \jmath = \alpha - \alpha'$, $\rho' v = \gamma - \gamma'$ and $\jmath' u + v \rho = \beta - \beta'$. And this means that the morphisms $(\alpha, \beta, \gamma)$ and  $(\alpha', \beta', \gamma')$ acting $\texttt z \To \texttt z'$ are \emph{homotopic} in the following sense:

\adef Two morphisms $f,g: C\To D$ between complexes are said to be homotopic when there is a morphism of complexes $u_i: C^i\To D^{i-1}$ so that
$$f_n - g_n = d^n u_n  + v_n d^n.$$
\zdef

The conclusion of of all this is that there is an homological notion, homotopic maps, that extends the equality notion in $\mathbb {SEX}$ to
$\mathbb {COM}$. We can test this notion: pick a projective (short) presentation $0\To \kappa \To \ell_1(\Gamma) \To X \To 0$ of the Banach space $X$. The identity of $X$ can be lifted in many ways to $\ell_1(\Gamma)$, and each lifting provides an arrow $(\alpha, \beta, 1)$. The point is that any two of these arrows  $(\alpha, \beta, 1)$ and $(\alpha', \beta', 1)$ are homotopic.  Thus, the identification of a space $X$ with its (short) projective presentations and operators with triples is faithful in the homotopic sense. Move ahead with this idea, do yourself a favour and delete the word ``short" in the last sentence (replacing triple by morphism of complexes) and show that whatever two projective presentations of a space $X$ are homotopic (or whatever two injective presentations, for that matter). And we are ready to conform the category $\mathbb {KOM}(\mathbf B)$ having complexes of Banach spaces as objects and homotopic equivalence classes of morphisms of complexes as morphisms.

\subsection{Exactness and how to measure it} Exactness is however important. A short exact sequence $0\To Y \stackrel{\jmath} \To Z \stackrel{\rho}\To X\To 0$ in $\mathbf{Ban}$ is, by the virtues of the open mapping theorem, a Banach space $Z$, a subspace $\jmath[Y]$ and the corresponding quotient $X=\rho[Z]= Z/\jmath[Y]$. A complex, even a short one, is... whatever it is. Homology attempts to measure how large the deviation of a complex from exactness is. Keep in mind that a complex could be exact at some point and not exact at another point. Fix a complex (we do not distinguish here between $i>0$ and $i<0$)
$$\begin{CD} \cdots @>>> C^{i+1} @>d^{i}>> C^{i} @>{d^{i-1}}>> C^{i-1}@>>> \cdots\end{CD} $$and define the $i^{th}$-homology group (vector space in our case)  $H^i(C)=\ker d^{i-1}/\;\mathrm{Im} \;d^i$. Of course that the complex $C$ is an exact sequence if and only if $H^i(C)=0$ for all $i$. When exactness occurs only at position $i$  one says that the complex $C$ is \emph{acyclic} at $i$ or exact at $i$. A morphism of complexes $f: C\To D$ generates morphisms $$H^if: H^i(C) \To  H^i(D)$$
in a natural way $ H^if (x+ \mathrm {Im}\; d^{i+1})= f_{i+1}(x) + \mathrm{Im} \;d^i$. It makes sense since $x\in \ker d^{i} \Rightarrow f_{i+1}(x)\in \ker d^{i}$ and since $f_{i+1}(x-y)= f_{i+1}d^{i+1}(z)= d^{i+1}f_{i+1}(z)$ then $x-y\in \mathrm {Im}\;d^{i+1}  \Rightarrow f_{i+1}x- f_{i+1}y\in
\mathrm {Im}\; d^{i+1}$. Now, different complex morphisms $f,g$ induce different maps between the homology groups. When two morphisms of complexes $f,g$ induce the same morphisms $H^if = H^ig$? When they are homotopic, which is not hard to show \cite[Lemma III.1.2]{guema}. Thus, a morphism of complexes $f:C\To D$ is homotopic to $0$ if and only if $H^i(f) =0$ for all $i$. The question of when a morphism $f$ induces a morphism $Hf$ so that all ${H^if}$ are isomorphisms is rather more elusive.

\adef A morphism of complexes $f$ is said to be a quasi-isomorphism when all ${H^if}$ are isomorphisms.\zdef

For instance, every morphism of complexes between exact sequences is a quasi-isomorphism. This does not mean that a morphism homotopic to $0$ has to be a quasi-isomorphism when the complexes are not exact!

Recapitulating, Banach space homology means working in a setting in which Banach spaces are identified with complexes ---long, large, possibly endless complexes---. These could be exact (or not), they could have been obtained using either projectives (complexes with $C_n=0$ for $n>0$) or injectives (complexes with $C_n=0$ for $n<0$) or not. An example of non-acyclic complex? This one: pick a Banach space $X$ and form
$$\begin{CD}  \dots 0@>>>0@>>> X @>>> 0@>>>0@>>>\dots\end{CD}$$
usually called $X[n]$ when $X$ is placed at position $-n$ (yes, there is a reason for that that we will encounter later). So, the seed idea is to to represent Banach spaces via complexes to then lift Banach functors ${\mathbf {Ban}}\To {\mathbf {Ban}}$ to functors ${\mathbf {Ban}}^{\mathbf {Ban}}\To {\mathbf {Ban}}^{\mathbf {Ban}}$, or at least to functors $\mathbb{COM}({\mathbf {Ban}})\To \mathbb{COM}({\mathbf {Ban}})$. If we try it, we will soon discover there is no natural way in which one can lift an operator $\tau:X\to X$ to a morphism between arbitrary exact complexes such as, say,
$$\begin{CD}
\dots@>>>0@>>> Y @ = Y@>>>0@>>>X @ = X@>>> 0\\
&&& &@V?VV@VV?V&   &@V{\tau}VV @VV{\tau}V\\
\dots@>>>0@>>> Z @ = Z@>>>0@>>>X @ = X@>>> 0
\end{CD}$$
let alone between complexes. However, projective presentations are more faithful in their representation of Banach spaces because of:

\begin{proposition} Every operator $\tau:X\To Y$ between two Banach spaces induces a morphism between whatever projective presentations of those spaces. Moreover, all morphism obtained that way are homotopic and quasi-isomorphisms. In particular, two projective presentations of $X$ are quasi-isomorphic.\end{proposition}

And how could one make a Banach functor $\mathcal F$ act between, say, projective presentations of two spaces $X,Y$? Well, it induces a functor
$\mathbb{COM}\To \mathbb{COM}$ acting term-by-term since it transforms complexes into complexes. Since it also transforms homotopic maps into homotopic maps, the induced functor also acts $\mathbb{KOM}\To \mathbb{KOM}$. A different thing is the behaviour of this induced map regarding quasi-isomorphisms. A functor $\mathcal F$ is said to be exact if it transforms short exact sequence into short exact sequences. It is clear that exact functors transform complexes quasi-isomorphic to $0$ (i.e. acyclic) into complexes quasi-isomorphic to $0$ (acyclic). Not as easy  \cite[Proposition III.6.2]{guema} and \cite[10.5.2]{weibel} is to show that:

\begin{proposition}\label{exactfunctor} An exact functor preserves quasi-isomorphisms.\end{proposition}

This suggests that another way to measure the deviation of a functor from exactness could be to check its behaviour with respect to quasi-isomorphisms. We will explore both ways in the next two sections: Section \ref{sec:derivation} explores how to measure the deviation from exactness; while
Section \ref{sec:derived} explores how to measure its behaviour on quasi-isomorphisms.

\subsection{Stop and hear the sound of waves: Homological dimension} The words still resonate in our minds: a long, large, possibly endless complex to represent a Banach space $X$. Yes, why not. With the only added problem that (almost) nothing, zero, zip, zilch, nada is known about even the simplest question: is actually infinite the projective presentation? Or the injective presentation for that matter. To the best of our knowledge, only Wodzicki \cite{wod} and recently Cabello, Castillo and Garc\'ia \cite{ccg} considered the topic. Let us define the  projective dimension $\mathfrak p(X)$ of a Banach space $X$ as the smallest $n$ for which $\ker d^n$ is projective; analogously,  the injective dimension $\mathfrak i(X)$ is is the smallest $n$ for which the cokernel space of $d^n$ is injective. It will take us no place to list all what is currently known:
\begin{itemize}
\item $\mathfrak p(X)=0$ if and only if $X$ is projective.
\item There is a quotient $\mathcal B$ of $c_0$ having $\mathfrak p(\mathcal B^*)=1$.
\item $\mathfrak p(X)=\infty$ if $X$ is an $\mathcal L_1$-space not projective.
\item If $\mathcal K$ denotes the Kadec space \cite{kade, pelcuni, p-w} then $\mathfrak p( \mathcal K)= \infty$.
\end{itemize}
As remarked in \cite{wod}, it is expected these dimensions to be $\infty$ for most ``classical'' spaces. Even if there is not a sniff of evidence (apart from many unsuccessful attempts) we conjecture that Hilbert spaces have infinite homological (projective or injective) dimension.

\section{Derivation}\label{sec:derivation}
Let us say it once again: exact sequences are important. Thus, exact functors are interesting. A good number of current-life Banach functors like \begin{itemize}
\item The duality functor $X\To X^*$, which is exact thanks to the Hahn-Banach theorem.
\item Functors of $p$-summable sequences: $X\To \ell_p(X)$ for $1\leq p\leq \infty$.
\item The ultraproduct functor.
\item The residual functor.
\end{itemize}
are exact, while several other important functors, such as
\begin{itemize}
\item Hom functors $\mathfrak L_X$ and $\mathfrak L^X$.
\item In general, the functors $\mathcal A(\cdot, X)$ and $\mathcal A(X, \cdot)$ for most operator ideals $\mathcal A$.
\item Tensor functors $\otimes_X$.
\item Functors of weakly $p$-summable sequences: $X\To \ell_p^w(X)$ for $1\leq p < \infty$.
\end{itemize}
are not. What occurs when a functor $\mathcal F$ is not exact? Well, Banach functors tend to be ``left-exact", in the sense that the image of a short exact sequence $0\to Y \to Z \to X\to 0$ is an exact sequence $0\to \mathcal F Y \to \mathcal F Z \to \mathcal F X$. What homology does in this case, or would like to do, is to obtain new functors $\mathcal F_1, \mathcal F_2 ,....$ that provide an exact sequence $$0\To \mathcal F Y \To \mathcal F Z \To \mathcal F X \To \mathcal F_1 Y \To \mathcal F_1 Z \To \mathcal F_1 X \To \mathcal F_2 Y \To \mathcal F_2 Z \To\mathcal F_2 X \To \cdots$$
in which the functor $\mathcal F_{n+1}$ measures how much $\mathcal F_n$ diverges from exactness, and for this reason $\mathcal F_{n+1}$ is called the derived functor of $\mathcal F_n$, and the process of passing from $\mathcal F$ to $\mathcal F_1$ is called derivation. How does one derive a functor? Not so simple: only \emph{certain} functors on \emph{certain} categories can be derived; and, according to MacLane \cite[XII, \S 9]{mac}: ``A standard method is: Take a resolution, apply a covariant functor, take homology of the resulting complex. This gives a connected sequence of functors, called the derived functors". Let's see this in action, with all due assumptions on a functor $\mathcal F$ to make sense of what follows: to obtain the left-derivation, pick an object $X$ in a suitable category, form a projective resolution
$$\begin{CD} \cdots @>>> \mathcal P_3 @>{d_3}>> \mathcal P_2 @>{d_2}>> \mathcal P_1 @>{d_1}>> X @>>> 0,\end{CD}$$
apply the functor $\mathcal F$ and get a \emph{complex} {\bf FX}
$$\begin{CD} \cdots @>>> \mathcal F\mathcal P_3 @>{\mathcal Fd_3}>> \mathcal F\mathcal P_2 @>{\mathcal Fd_2}>>\mathcal F\mathcal P_1 @>{\mathcal Fd_1}>> \mathcal FX @>>> 0\end{CD}$$
 The homology group $H^n({\bf FX})$ can be used to define the $n^{th}$ (left) derived functor $L_n\mathcal F$ of $\mathcal F$ at $X$. $L_1\mathcal F$ is a functor because arrows $f: A\To B$ yield a morphism $ H^n({\bf FA}) \To  H^n({\bf FB})$ well defined since any two liftings of $f$ are homotopic. A similar process worked out with injective presentations yields the right derived functors. There are other
 contructions that can be called ``derived functor". The existence of a long homology exact sequence formed by derived functors still depends on the existence of connecting morphisms $\mathcal F_nX \To \mathcal F_{n+1} Y$. Let us show how they appear in the construction above of left derivation: given an exact sequence $0\to Y \to Z \to X \to 0$, construct connected projective resolutions of the spaces $Y,Z, X$ forming a commutative diagram
$$
\xymatrixrowsep{1.5pc}
\xymatrix{
& \cdots\ar[d] & \cdots\ar[d]  & \cdots\ar[d]  &\\
0\ar[r]& \mathcal P_2^Y\ar[r]\ar[d] & \mathcal P_2^Z\ar[r]\ar[d] & \mathcal P_2^X\ar[r]\ar[d]& 0\\
0\ar[r]& \mathcal P_1^Y\ar[r]\ar[d] & \mathcal P_1^Z\ar[r]\ar[d] & \mathcal P_1^X\ar[r]\ar[d]& 0\\
0\ar[r]& Y \ar[r]\ar[d] & Z\ar[r]\ar[d] & X\ar[r]\ar[d]& 0\\
& 0 & 0 & 0&
}$$
This yields an exact sequence of complexes $0\To {\bf Y}\To {\bf Z}\To {\bf X}
\To 0$. Now, the necessary screw that supports the construction is provided by the observation that each map $d_n$ in a complex {\bf C} induces a map $\overline d_n: \coker d_{n+1}\to \ker d_{n-1}$ so that $\ker \overline d_n = H^n({\bf C})$ and $\coker \overline d_n = H^{n-1}({\bf C})$ as it is easy to see \cite[IV Lemma 2.2.]{hiltstam}. One thus gets a commutative diagram
$$\xymatrix{
&H^n\ar[r]({\bf Y})\ar[d]& H^n({\bf Z})\ar[r]\ar[d]& H^n({\bf X})\ar[d]\\
&\coker d_{n+1}\ar[r]\ar[d]_{\overline d_n}& \coker d_{n+1}\ar[r]\ar[d]_{\overline d_n} & \coker d_{n+1}\ar[r]\ar[d]^{\overline d_n} &0\\
0\ar[r]&\ker d_{n-1}\ar[r]\ar[d]& \ker d_{n-1}\ar[r]\ar[d]& \ker d_{n-1}\ar[d]\\
&H^{n-1}({\bf Y})\ar[r]& H^{n-1}({\bf Z})\ar[r] & H^{n-1}({\bf X})}$$
The Snake lemma yields a connecting morphism $H^n({\bf X})\To H^{n-1}({\bf Y})$. To conclude, when projective and injective presentations are available, right and left derivations are naturally equivalent, and thus the derived functor is unique. In $p$-Banach spaces, at least projective derivation is possible. But that is no longer the case of quasi Banach spaces. So, a nice question, to the best of my knowledge posed by F\'elix Cabello, is: how to derive quasi Banach functors?

\subsection{Pedestrian, likely misleading, examples.} Let's pick as template the case of one of the basic Banach functors: $\mathfrak L_A$. To understand its derivation, we follow HMBST: let us make a strategic move and attempt to guess how, given an exact sequence $0\To Y  \To Z \To X \To 0$, we could continue the sequence
$$\begin{CD} 0@>>> \mathfrak L (A, Y) @>>> \mathfrak L (A, Z) @>>> \mathfrak L (A, X) \end{CD}$$
keeping exactness. This entails finding a vector space $V$ and a linear map $L: \frak L(A,X)\To V$ whose kernel is the space of  operators $\frak L(A,X)$ that lift to $Z$. The simplest solution could be: $L$ should send $\tau\in \frak L(A, X)$ to the lower exact sequence in the pullback diagram
$$
\xymatrix{
0 \ar[r] & Y\ar[r] \ar@{=}[d] & Z \ar[r] &X\ar[r] & 0\\
0 \ar[r] & Y\ar[r] & \PB \ar[r] \ar[u] &A\ar[r] \ar[u]_\tau & 0}
$$
because $L(\tau)=0$ precisely when $\tau$ can be lifted to an operator $\frak L(A, Z)$. According with this, our best candidate for $V$ is the space $\Ext(A, Y)$ of exact sequences $0\To Y  \To \diamondsuit \To A \To 0$, modulo the standard equivalence relation. The manifest implication is that the continuation can be provided by the functor $X \mapsto \Ext(A, X)$. In other words, it is reasonable to think that $\Ext(A, \cdot)$ is \emph{a} derived functor of $\mathfrak L(A, \cdot)$ and $\Ext(\cdot, A)$ is \emph{a} derived functor of $\mathfrak L(\cdot, A)$. The point to beware with is that all of this we have done works exactly the same in quasi Banach spaces, even if no projective or injective elements exist in that category and thus there is no clear justification for the assertion ``$\Ext$ is \emph{the} derived functor of $\mathfrak L$". The way of looking at Kan extensions of $\Ext_{\mathbf{Ban}}$ as reasonable candidates as derivations of $\mathfrak L$ in $\mathbf {QBan}$ was what we did, without a stratospheric success, in Section \ref{forwe}.

Let us sum up known facts
 \begin{itemize}
\item The derived functor of an Exact functor should be $0$.
\item $\Ext_{\mathbf{Ban}}$ is the derived functor of $\mathfrak L: \mathbf{Ban}\To \mathbf { Vect}$.
\item The derived functor of $\otimes$ is called Tor.
\item The derivation of operator ideal functors $\mathcal A$ opens the topic of relative homology \cite{medit}
\item Derivation continues, and the derived functor of $\Ext_{\mathbf{Ban}}$ is called $\Ext_{\mathbf{Ban}}^2$. In general, the derived functor of $\Ext_{\mathbf{Ban}}^n$ is denoted $\Ext_{\mathbf{Ban}}^{n+1}$. The nature and properties of $\Ext_{\mathbf{Ban}}^n$ are, by the times of writing these lines, a good mystery.
\item The paper \cite{ccg} was written to uncover the Yoneda approach to $\Ext^n$ in Banach and, to some extent, in quasi Banach spaces.
\item The central result in the study of exact sequences and twisted sums of Banach spaces is that $\Ext(\ell_2, \ell_2)\neq 0$ \cite{elp, kaltpeck}. In \cite{ext2} it is shown that also $\Ext^2(\ell_2,\ell_2)\neq 0$. We conjecture that $\Ext^n(\ell_2, \ell_2)\neq 0$ as well.
\end{itemize}

\section{A necessary, maybe, step: Localization}

One feels better off working in a Banach space theory in which all one dimensional spaces are the same. This has effects even at the mathematics level since the set of all real one dimensional Banach spaces up to isometries has just one element $\R$ while the set of all real one dimensional Banach spaces is not even a set (since it essentially is the whole mathematical universe). After that one has grown accustomed to believe that all isometric Banach spaces are the same. Which, in combination with the outstanding result that all Hilbert spaces of the same dimension are isometric sooner or later makes one slip into the idea that all Hilbert spaces are the same. Let us continue from success to success until the final defeat: who has not dreamed of working in a Banach space theory in which all reflexive spaces are the same? Or all finite dimensional spaces are the same? From an operator ideal point of view, we could reformulate this as: pick an operator ideal $\mathcal A$ and create a Banach space category in which all elements of $\mathcal A$ are isomorphisms. Could we do that?
Of course not, but for a reason intrinsic to operator ideals: $0$ is always an element of the ideal. Thus, operator ideals are bad for localization, which means that other operator structures are needed. A similar phenomenon was pointed at in \cite{3sp}: operator ideals are bad for ``three-operator problems"; in this case, are the semigroups of operators which seem to be adequate. \smallskip

Fix a category $\mathbf{B}$ and a family $\mathcal A$ of morphisms of $\mathbf{B}$. The localized category $\mathbf{B}[\mathcal A^{-1}]$ has as objects the same as $\mathbf{B}$, while morphisms of $\mathbf{B}[\mathcal A^{-1}]$ are equivalence classes of arrows according to the following rules \cite[1.1]{gazi}:
\begin{itemize}
\item Introduce for each $a: X\to Y$ in  $\mathcal A$ an arrow $a^{-1}: Y\to X$ in $\mathbf{B}[\mathcal A^{-1}]$.
\item Each arrow $f:X\to Y$ in $\mathbf{B}$ is an arrow $f:X\to Y$ in $\mathbf{B}[\mathcal A^{-1}]$.
\item A finite sequence of arrows in $\mathbf{B}[\mathcal A^{-1}]$ each of them ending where the previous one ends is an arrow in $\mathbf{B}[\mathcal A^{-1}]$.
\item The arrows $a^{-1} a$ and $a a^{-1}$ are equivalent to the identity in $\mathbf{B}[\mathcal A^{-1}]$.
\item Two consecutive arrows is equivalent to their composition.
\end{itemize}

The localization $\mathbf{B}[\mathcal A^{-1}]$ has an associated functor $\mathcal Q: \mathbf{B} \To \mathbf{B}[\mathcal A^{-1}]$ sending an arrow $f:X\to Y$ of $\mathbf{B}$ to the morphisms $f:X\to Y$ of $\mathbf{B}[\mathcal A^{-1}]$. In this way, (the class of) $a:X\to Y$ in $\mathcal A$ is an isomorphism in $\mathbf{B}[\mathcal A^{-1}]$ with inverse (the class of) $a^{-1}:Y\to X$. The localized category comes accompanied with the corresponding universal property:

\begin{proposition} Let $\CAT C$ be a category and let $Q: \mathbf{B} \To \CAT C$ be a functor transforming morphisms in $\mathcal A$ into isomorphisms. Then $Q$ factorizes through $L_{\mathcal A}$ in the form
$$\xymatrix{
\mathbf{B} \ar[rr]^{\mathcal Q} \ar[d]_Q&& \mathbf{B}[\mathcal A^{-1}]\ar[dll]^F\\
\CAT C & }
$$
\end{proposition}
\begin{proof} The functor $F$ is cleanly defined as $F(X)=X$ for objects, $F(f)=f$ if $f$ is an arrow of $\mathbf{B}$ and
$F(a^{-1})=Q(a)^{-1}$.\end{proof}

Suitably chosen classes $\mathcal A$ allow nicer descriptions of the localized categories:

\adef A class $\mathcal A$ of morphisms will be called (left) \emph{localizing} if it enjoys the following stability properties:
\begin{itemize}
\item Is closed under composition.
\item $1_X\in \mathcal A$ for every object $X$.
\item In a pullback diagram
$$\xymatrix{
A\ar[r]^a& B\\
\PB\ar[r]_{\underline a}\ar[u]& C\ar[u]}$$
if $a\in\mathcal A$ then also $\underline a \in \mathcal A$.
\item Given two arrows $\tau, \eta: X\to Y$ there exists $a\in \mathcal A$ such that $a\tau=a\eta$ if and only if there exists $a\in \mathcal A$ such that $\tau a = \eta a$.
\end{itemize}\zdef

Now, a localizing class $\mathcal A$ allows a rather visual approach to the localized category $\mathbf{B}[\mathcal A^{-1}]$: objects are the same as in $\mathbf{B}$ while a morphism $X\to Y$ in  $\mathbf{B}[\mathcal A^{-1}]$ is an equivalence class of roofs $(f,a)$ with $f$ a morphism and $a\in \mathcal A$:
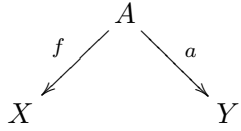
$$\xymatrix{&A\ar[dl]_{f}\ar[dr]^a&\\
X& &Y}$$
which corresponds to $\mathcal Q(f)\mathcal Q(a)^{-1}$. (Astute readers won't get fooled: if one imagines composition of $f$ with $a$ then $a$ should be on the left and $f$ on the right, so that
when $f=a$ ten $(a,a)$ is the identity of $X$ (not ``the opposite of the identity" as it occurs with our drawing). Since aesthetically is nicer to have $a$ on the right, my advice is simple: stop imagining). Two roofs $(f,a)$ and $(f',a')$ are equivalent if and only id there exists another roof $(f'', a'')$ making a commutative diagram
$$
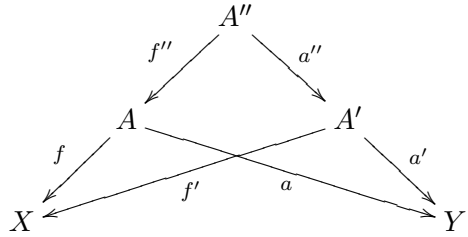
\xymatrix{&&A''\ar[dl]_{f''}\ar[dr]^{a''}&&\\
&A\ar[dl]_{f}\ar[drrr]_a&&A'\ar[dr]^{a'}\ar[dlll]^{f'}\\
X&&&&Y
}$$
The composition of two roofs $(f,a)$ and $(f',a')$ is the roof
$$
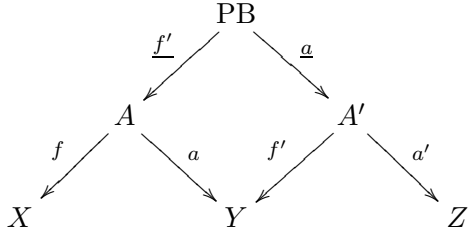
\xymatrix{&&\PB\ar[dl]_{\underline{f'}}\ar[dr]^{\underline{a}}&&\\
&A\ar[dl]_{f}\ar[dr]^a&&A'\ar[dr]^{a'}\ar[dl]_{f'}\\
X&&Y&&Z
}$$

The equivalence of roofs is actually an equivalence relation and composition respects this equivalence. To prove that is uncomplicated, but not exactly simple. Let us call $\widehat {\mathbf{B}}$ the category of roofs:\\

\textbf{Claim.} $\widehat {\mathbf{B}} = \mathbf{B}[\mathcal A^{-1}]$.\\

To prove it, it is enough to check that $\widehat {\mathbf{B}}$ endowed with the functor $\wedge: \mathbf{B} \To \widehat {\mathbf{B}}$ given by
$\wedge(X)=X$ for objects and $\wedge (f) = (f, 1_X)$ enjoys the same universal property: given a category $\CAT C$ and a functor $Q: \mathbf{B} \To \CAT C$ transforming morphisms in $\mathcal A$ into isomorphisms then $Q$ factorizes through $\wedge$ in the form
$$
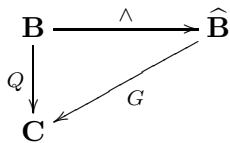
\xymatrix{
\mathbf{B} \ar[rr]^{\wedge} \ar[d]_Q&& \widehat {\mathbf{B}}\ar[dll]^G\\
\CAT C & }
$$
To construct $G$ we set $G(X)=X$ for objects and $G(f,a)= Q(f)Q(a)^{-1}$ (check that it is well defined) which clearly makes $G\wedge = Q$.\qed

We could have worked with the dual construction just working in the opposite category:

\adef A class $\mathcal A$ of morphisms will be called (right) \emph{localizing} if it enjoys the following stability properties:
\begin{itemize}
\item Is closed under composition.
\item $1_X\in \mathcal A$ for every object $X$.
\item In a pushout diagram
$$\xymatrix{
A\ar[r]^a\ar[d]& B\ar[d]\\
C\ar[r]_{\overline a}& PO}$$
if $a\in\mathcal A$ then also $\overline a\in \mathcal A$.
\item Given two arrows $\tau, \eta: X\to Y$ there exists $a\in \mathcal A$ such that $a\tau=a\eta$ if and only if there exists $a\in \mathcal A$ such that $\tau a = \eta a$.
\end{itemize}\zdef

This yields the (right) localized category $\mathbf{B}[\mathcal A^{-1}$]: objects are the same as in $\mathbf{B}$ while a morphism $X\to Y$ in  $\mathbf{B}[\mathcal A^{-1}$] is an equivalence class of floors$(f,a)$:
$$\xymatrix{X\ar[dr]_{f}& &Y\ar[dl]^a\\
&A&}$$
All the rest goes as it should. Another possible approach is to pick as $\mathcal A$ a class of what is called weak equivalences:
\adef A class $\mathcal A$ of weak equivalences is a class of morphisms in a category $\CAT C$ that contains the identity and satisfies the 2-out-of-6 property: whenever three arrows
$$\begin{CD} A@>u>> B @>v>> C  @>w>> D\end{CD}$$
are such that $wv$ and $vu$ belong to $\mathcal A$ then $u,v,w$ and $wvu$ belong to $\mathcal A$.
\zdef

A category $\CAT C$ equipped with a class of weak equivalences is usually called a \emph{homotopical category} and yields after localization a category $\CAT{HC}$ called the \emph{homotopy category}. There is an obvious localization functor $H_{\CAT C}: \CAT{C}\To \CAT{HC}$, and it has the universal property with respect to categories $\CAT D$ endowed with a functor $\eta: \CAT C \To \CAT D$ that transforms weak equivalences into isomorphisms: $\eta$ factorizes through $H_{\CAT C}$.

Of course, the identification of right or left localizing classes of operators in Banach spaces is an untouched topic. On the other hand, isomorphisms form quite obviously a class of weak equivalences. When working in categories of complexes $\mathbb{KOM}$ then quasi-isomorphisms usually conform the class of weak equivalences.

\section{Derived category}\label{sec:derived}

The notion of derived category $\mathbb {D}(\CAT C)$ of a given a suitable category $\CAT C$ is simple to describe but not as simple to handle:  it is the localization of $\mathbb{KOM}(\CAT C)$ with respect to quasi-isomorphisms. One could also have said that $\mathbb {D}(\CAT C)$ is the localization of $\mathbb{COM}(\CAT C)$ with respect to quasi-isomorphisms: if $\mathcal Q$ is the localization functor then whenever $f$ is homotopic to $g$ then $\mathcal Q(f)= \mathcal Q(g)$. However the problem is not that, the problem is that the notion is not so easy to swallow, is it?\\

Be as it may, where we are pointing at is: does a functor $\mathcal F: \CAT C\To \CAT C$ induce a functor $\mathbb {D}(\CAT C)\To \mathbb {D}(\CAT C)$? We see that $\mathcal F$ induces, by plain term-by-term application to complexes, a functor $\mathbb {C}(\mathcal F): \mathbb {COM}(\CAT C)\To \mathbb {COM}(\CAT C)$. This functor respects homotopy and thus it induces a functor $\mathbb {K}(\mathcal F): \mathbb {KOM}(\CAT C)\To \mathbb {KOM}(\CAT C)$. Now, given an object $C$, represented by a complex $\mathscr C$ if we apply $\mathcal F$ componentwise we obtain a complex $\mathcal F \mathscr C$ that represents $\mathcal F C$. If the complex $\mathscr C$ is acyclic and $\mathcal F$ is exact then $\mathcal F \mathscr C$ is acyclic too. Recall from Proposition \ref{exactfunctor} that
 exact functors transform quasi-isomorphism into quasi-isomorphisms, and thus it induces a functor $\mathbb{D}(\mathcal F): \mathbb{D}(\CAT C)\To \mathbb{D}(\CAT C)$. In the general case, when $\mathcal F$ is not exact, then an induced functor $\mathbb {D}(\CAT C)\To \mathbb {D}(\CAT C)$ could exist or not. But, even when it does, it is not necessarily exact and, therefore, the homology groups $H^i(\mathbb D(\mathcal F)(\textbf C))$ are not necessarily $0$. The idea then is that  $H^i(\mathbb D(\mathcal F)(\textbf C))$ is the $i^{th}$-derived functor of $\mathcal F$ at $C$. After all, if $\mathbb D(\mathcal F)$ respects quasi-isomorphisms, these homology groups are isomorphic no matter which projective presentation of $C$ has been chosen. According to this, the derived functors of an exact functor $\mathcal F$ are $0$, as they should. The return journey would then be that the derived functors of a functor $\mathcal F$ (in the standard sense earlier described) conform a suitable extension of $\mathcal F$ to $\mathbb D(\CAT C)$. Of course that we have not arrived thus far to think that ``a suitable extension" could be anything different from a Kan extension.\medskip

To see these ideas in action the problem is, of course, the nature of the derived category which, in turn, depends on the nature of the starting category. That is what makes so complicated the work of homologists. But we are just interested in what occurs with $\mathbf{Ban}$ and $\mathbf {QBan}$ so, it is enough for us to outdo the Night King instead of worrying about what happens with the iron throne. Some of the byways one often follows in order to have a manageable representation of $\mathbb D(\CAT C)$ are:

\begin{itemize}
\item Instead of using \emph{all} complexes, work only with complexes on the positive or negative range. This leads to $\mathbb{COM}^{+}$ (complexes with $C_i=0$ for all $i<0$) and then $\mathbb{KOM}^{+}$. Analogously with  $\mathbb{COM}^{-}$ (complexes with $C_i=0$ for all $i>0$) and then  $\mathbb{KOM}^{-}$. The functors induced by some $\mathcal F$ will be then either $\mathbb C^+(\mathcal F)$ or $\mathbb C^-(\mathcal F)$ and $\mathbb K^+(\mathcal F)$ or $\mathbb K^+(\mathcal F)$ as well.
\item That produces also the derived categories $\mathbb{D}^{-}$ and $\mathbb{D}^{+}$ .
\item Instead of using arbitrary complexes, use only a limited amount of them. For instance, complexes formed only by injective (or projective) objects (not necessarily acyclic!). This produces the categories $\mathbb{KOM}(\CAT C, inj)$ or $\mathbb{KOM}(\CAT C, proj)$.
\item Using only injective or projective presentations yields a certain disruption: all morphism of complexes between exact complexes are quasi-isomorphisms.
\end{itemize}

Keeping these strategies in mind one can obtain:

\begin{enumerate}
\item \cite[Theorem III.4.4]{guema} The class of quasi-isomorphisms is localizing in $\mathbb{KOM}(\CAT C)$. This means that the approach to $\mathbb D(\CAT C)$ via roofs is possible. It is also localizing in $\mathbb{KOM}^+(\CAT C, inj)$
\item \cite[Theorem 10.4.8]{weibel} When $\CAT C$ has enough projectives then $\mathbb D^{-}(\CAT C)$ exists and is equivalent to the (full) subcategory $\mathbb{KOM}^{-}(\CAT C, proj)$.
\item \cite[Theorem 10.4.8]{weibel} \cite[Theorem III.5.21]{guema} When $\CAT C$ has enough injectives then $\mathbb D^{+}(\CAT C)$ exists and is equivalent to the (full) subcategory $\mathbb{KOM}^{+}(\CAT C, inj)$.
\item (please go to Section \ref{heart} to discuss about what ``Abelian category" means) An Abelian category $\CAT C$ can be seen as the full subcategory of $\mathbb D(\CAT C)$ formed by the complexes $\mathbf C$ with $H^i(\mathbf C)=0$ for $i\neq 0$.
\end{enumerate}

The ``exists" assertion coda in (2) and (3) may surprise devil-may-care readers. It means ``it exists in our universe": when $\CAT C$ is a small category, localizations exist. Otherwise ... See \cite[Remark 10.3.3]{weibel} for a discussion on this topic when $\CAT C$ is not small. Once again, we have not arrived thus far to pretend believing that $\mathbb{D}(\CAT C)$ is not characterized by some universal property:
\begin{enumerate}
\item[(4)] \cite[Definition - Theorem II.2.1]{guema} If $\CAT C$ is an Abelian category there exists a category $\mathbb D(\CAT C)$ and a functor $\mathcal Q: \mathbb{KOM}(\CAT C) \To \mathbb{D}(\CAT C)$ with the following properties
\begin{itemize}
\item $\mathcal Q$ transforms quasi-isomorphisms into isomorphisms.
\item Any functor $\mathcal F: \mathbb{KOM}(\CAT C) \To {\CAT D}$ transforming quasi-isomorphisms into isomorphisms factors through $\mathcal Q$.
\end{itemize}
\end{enumerate}
This universal property is at the root of assertions (2) and (3) above: if one is able to work only with some type of exact complexes in such a way that quasi-isomorphisms become isomorphisms then the corresponding $\mathbb{KOM}$ is the derived category. But, as many times in life, in categorical life whatever is not too much is not even enough.

\section{Derived functor} Different things, in different context, are called derived functors. In the context of complexes and the derived category one has \cite[III.5.6]{guema}:

\begin{theorem} Let $\mathcal F: \CAT C \To \CAT D$ be an additive functor acting between Abelian categories.
\begin{itemize}
\item The right derived functor $R\mathcal F$ of a left exact functor $\mathcal F$ is the left Kan extension of $\mathcal Q  \; \mathbb K^+(\mathcal F): \mathbb{KOM}^+(\CAT C)\To \mathbb D^+(\CAT D)$ through $\mathcal Q$. This means that $R\mathcal F$ makes a commutative square
$$\xymatrix{\mathbb {KOM}^+(\CAT C) \ar[ddd]_{\mathbb {K}^+(\mathcal F)} \ar[rrr]^{\mathcal Q}& &&\mathbb {D}^+(\CAT D)\ar[ddd]^{R\mathcal F}\\
&&\\
&\;\ar[ur]^\epsilon&\\
\mathbb {KOM}^+(\CAT C) \ar[rrr]_{\mathcal Q} &&& \mathbb {D}^+(\CAT D)}$$and there is a natural transformation $\epsilon: {\mathcal Q}\; \mathbb {K}^+(\mathcal F)\To R\mathcal F \circ\mathcal Q$ yielding the uniqueness property with respect to that diagram.

\item The left derived functor of a right exact functor $\mathcal F$ is the right Kan extension of $\mathcal Q \;\mathbb K^-(\mathcal F):
\mathbb{KOM}^-(\CAT C)\To \mathbb D^-(\CAT D)$ through $\mathcal Q$. This means that $L\mathcal F$ makes a commutative square
$$\xymatrix{\mathbb {KOM}^-(\CAT C) \ar[ddd]_{\mathbb {K}^-(\mathcal F)} \ar[rrr]^{\mathcal Q}&& &\mathbb {D}^-(\CAT D)\ar[ddd]^{L\mathcal F}\\
&&\;\ar[dl]^\eta&\\
&&\\
\mathbb {KOM}^-(\CAT C) \ar[rrr]_{\mathcal Q} &&& \mathbb {D}^-(\CAT D)}$$and there is a natural transformation $\eta: L\mathcal F \;\mathcal Q \To {\mathcal Q}\;\mathbb {K}^-(\mathcal F)$ yielding the uniqueness property with respect to that diagram.
\end{itemize}
\end{theorem}

In the context of homotopic categories (i.e., categories $\CAT C$ endowed with weak equivalences so that one can construct the category $\CAT{HC}$) the formulation is similar \cite[6.4]{emily}:
\adef Let $\mathcal F: \CAT C \To \CAT D$ be a functor between homotopical categories
\begin{itemize}
\item The Right Kan extension of $H_{\CAT D} \mathcal F$ through $H_{\CAT C}$ defines the so-called \textbf{total} left derived functor of $\mathcal F$
\item The Left Kan extension of $H_{\CAT D} \mathcal F$ through $H_{\CAT C}$ defines the so-called \textbf{total} right derived functor of $\mathcal F$
\end{itemize}\zdef\medskip

All this, however, leaves us with a bittersweet taste since our expectations were that derivation of a functor $\mathcal F: \CAT C\To \CAT D$ would be a functor $\CAT C\To \CAT D$. Let us see two ways to proceed when the target category $\CAT D$ is Abelian, so that taking homology ``groups" provides objects of $\CAT D$.\medskip

\begin{enumerate}
\item Define as the $n^{th}$-right derived functor of $\mathcal F: \CAT C\To \CAT D$ the composition
$$
\begin{CD}
\CAT C @>[n]>> \mathbb{KOM}(\CAT C)@>{R\mathcal F}>> \mathbb{KOM}(\CAT D)@>{H^0}>>  \CAT D\end{CD}$$
where
\begin{itemize}
\item $R\mathcal F$ is the right derived functor as described above.
\item $[n]$ is the functor $X\To X[n]$ described in Section \ref{homo} in which t $X[n]$ is the non-acyclic complex $ \dots \To  0\To  0\To  X\To 0 \To 0 \To \dots$ with $X$ placed at position $-n$.
\item $H^0$ is just taking the homology group at position $0$.
\end{itemize}
and define the $n^{th}$-left derived functor of $\mathcal F: \CAT C\To \CAT D$ will be defined as the composition
$$
\begin{CD}
\CAT C @>[n]>> \mathbb{KOM}(\CAT C)@>{L\mathcal F}>> \mathbb{KOM}(\CAT D)@>{H^0}>>  \CAT D\end{CD}$$
where
\begin{itemize}
\item $L\mathcal F$ is the right derived functor as described above.
\item $[n]$ is the functor $X\To X[n]$ described in Section \ref{homo} in which t $X[n]$ is the non-acyclic complex $ \dots \To  0\To  0\To  X\To 0 \To 0 \To \dots$ with $X$ placed at position $-n$.
\item $H^0$ is just taking the homology group at position $0$.
\end{itemize}\medskip

\item Given a right exact functor between suitable categories admitting enough projectives, a $n^{th}$-left derived functor would be provided by the composition
$$
\begin{CD}
\CAT C @>>> \mathbb{KOM}^-(\CAT C,\text{proj})@>{\mathbb{K}^- (\mathcal F)}>> \mathbb{KOM}^-(\CAT D, \text{proj})@>{H^n}>>  \CAT D\end{CD}$$
while a $n^{th}$-right derived of a left exact functor between suitable categories admitting enough injectives would be provided by the composition
$$
\begin{CD}
\CAT C @>>> \mathbb{KOM}^+(\CAT C, \text{inj})@>{\mathbb{K}^+ (\mathcal F)}>> \mathbb{KOM}^+(\CAT D, \text{inj})@>{H^n}>>  \CAT D\end{CD}$$
\end{enumerate}

In those cases, the natural inclusions $\CAT C \To \mathbb{KOM}^-(\CAT C, \text{proj})$ and $\CAT C \To \mathbb{KOM}^+(\CAT C, \text{inj})$ are the obvious ones. In general, a category $\CAT C$ can be seen in different ways inside $\mathbb{KOM}(\mathbf {\CAT C})$, but the most natural is to represent is as the full subcategory of $H^0$-complexes; i.e., complexes $\mathscr C$ for which $H^i(\mathscr C)=0$ for all $i\neq 0$. Further composition with the localization functor $\mathcal Q: \mathbb{KOM} \To \mathbb D(\CAT C)$ yields a full subcategory equivalent to $\CAT C$.

\section{Back to Ban}

 How do all these ideas work in Banach spaces? Because it is always the underlying problem that $\mathbf {Ban}$ is not an Abelian category (Section \ref{heart}). Since it has enough injectives and projectives, we can well use injectives and work with $\mathbb D^{+}(\mathbf {Ban}) = \mathbb{KOM}^{+}(\mathbf {Ban})= \mathbb{KOM}^{+}(\mathbf {Ban}, inj)$ or use projectives and  work with $\mathbb D^{-}(\mathbf {Ban}) = \mathbb{KOM}^{-}(\mathbf {Ban})= \mathbb{KOM}^{-}(\mathbf {Ban}, proj)$. All this is fairly standard. Or we could construct the derived category, which is not standard. Let us present first the classical approaches, with a twist, to then exhibit the construction via the derived category. We will do all those things keeping the functor $\mathfrak L_X$ in focus.\medskip

In all cases, a warning must be made right now: the ``extraction" functors $H^0, H^n$ do not take values in $\mathbf{Ban}$. We would very much like that were true but it is not. Guess why? Yes, because $\mathbf{Ban}$ is not an Abelian category and thus, as a rule, when one has an operator $T: X\To Y$, the quotient $Y/ T[X]$ is not a Banach space. This generates a crisis but, according with the original meaning of the word, gives us a non insignificant opportunity: homology groups (vector spaces) are not Banach spaces, but they are semi-Banach spaces; namely, everything is Ok except that we must work with a seminormed not-necessarily-Hausdorff topology. Therefore, if we call $\mathbf{sBan}$ the category of semi Banach spaces, the functors $H^n$ act between whatever category of complexes one is using and $\mathbf{sBan}$.

\subsection{The standard approach, with a twist}

The standard left derivation procedure of $\mathfrak L^E$ is as follows: given the Banach space $X$, pick a projective presentation $0\to \kappa\to P\to X\to 0$, then a projective presentation of $\kappa$, all of which we can assemble as $\kappa'\to P' \to \kappa\to P\to X\to 0$
and this forms an exact sequence $\kappa'\To P' \To  P\To X\to 0$ that we will call $\mathscr P$ (we do not care what comes next to get the first derivation). Apply $\mathfrak L^E$ to get $\mathfrak L^E(\mathscr P)$, namely $0\to \mathfrak L(X, E) \To \mathfrak L(P, E) \To \mathfrak L(\kappa, E) \To \mathfrak L(P', E)$ that we reassemble as the complex $0\To \mathfrak L(X, E) \to \mathfrak L(P, E) \To \mathfrak L(P', E)$
that will be called $\mathfrak L^E(\mathscr P)$. Therefore $L(\mathfrak L_E)(X) = H^1(\mathfrak L^E(\mathscr P))=\mathfrak L(\kappa, E)/\sim$ , where $\sim$ is the equivalence relation $\tau\sim \eta \Longleftrightarrow $$ \tau-\eta$ can be extended to an operator $P\to E$. The reader can take some  time to check that what has been done does not depend on the projective presentations (or read it in HMBST Chapter 2). Thus, the basic results about the pushout construction imply that precisely $\mathfrak L(\kappa, E)/\sim \;= \Ext_{\mathbf{Ban}}(X, E)$. This is the way we are certain that
$$\mathfrak L(\cdot, E)' = \Ext_{\mathbf{Ban}}(\cdot, E)$$

The construction above has a weak point: what has been actually shown is that $\mathfrak L(\cdot, E)' = \mathfrak L(\kappa, E)/\sim$, which is a bad formula because $\kappa$ depends on the projective space $P$ and the quotient map chosen. So, we would be much more happy if we could write $\mathfrak L(\cdot, E)' = \mathfrak L(\kappa(X), E)/\sim$, where now $\kappa$ is some Banach functor. That can be done using the functor $\co_z(\cdot)$, different forms of whose construction can be found in \cite{moretesis,morext,sob,cacas,duality}, and in HMBST, Section 3.10. The key point is that there is a Banach functor $co_z(\cdot)$ such that for each Banach space there is an exact sequence $\diamondsuit$ with the form $0\To co_z(X)\To \diamondsuit\To X\To 0$ such that every element $\spadesuit$ in $\Ext_{\mathbf{Ban}}(X, E)$ is equivalent to a pushout of $\diamondsuit$. In other words, there is a commutative diagram
$$\xymatrix{
0\ar[r]& co_z(X)\ar[r]\ar[d]_{\phi_{\spadesuit}}&\diamondsuit \ar[d]\ar[r]& X\ar[r]& 0\\
0\ar[r]& E\ar[r]& \spadesuit \ar[r]& X\ar@{=}[u]\ar[r]& 0}$$
In this way, again by the general properties of the pushout construction, the correspondence $\spadesuit \To \phi_{\spadesuit}$ establishes an identification
$$\Ext_{\mathbf{Ban}}(X, E) =  \mathfrak L(co_z(X), E)/\sim$$

This construction can be iterated to provide

$$\Ext_{\mathbf{Ban}}^2(X, E) = \Ext_{\mathbf{Ban}}(co_z(X), E) =  \mathfrak L(co_z(co_z(X)), E)/\sim$$

and analogously for higher order derived $\Ext^n$ spaces. See \cite{ccg} for details. A kind of dual construction works in the same way for injective presentations and right derivation.

\subsection{Using the derived category}

In the description above one gets$$
\begin{CD}
\mathbf{Ban} @>[n]>> \mathbb{KOM}(\mathbf{Ban})@>{R\mathfrak L_X}>> \mathbb{KOM}(\mathbf{Ban})@>{H_0}>>  \clubsuit\end{CD}$$
where we would like to place $\mathbf{ Ban}$ instead of $\clubsuit$, but we can't.
$$\begin{CD}
\mathbf{Ban} @>[n]>> \mathbb{KOM}(\mathbf{Ban})@>{R\mathfrak L_X}>> \mathbb{KOM}(\mathbf{Ban})@>{H_0}>>  \mathbf{sBan}\end{CD}$$

and thus what we can obtain is that if $\mathfrak L_X^{(n)} = H^0R\mathfrak L_X[n]$ is the $n^{th}$-right derived functor of $\mathfrak L_X$ above then $\mathfrak  L_X^{(n)} (Y)$ is a semi-Banach space. But, who is this $\mathfrak  L_X^{(n)} (Y)$ ?

\begin{proposition} $\Ext_{\mathbf {Ban}}^n(X,Y) = \Hom_{\mathbb{KOM}(\mathbf{Ban})}(X[0], Y[n])$
\end{proposition}

Instead of making a really annoying proof, let us explain why $\Ext_{\mathbf {Ban}}(X,Y) = \Hom_{\mathbb{KOM}(\mathbf{Ban})}(X[0], Y[1])$: because exact sequences $0\To Y \to Z\to X\To 0$ correspond to morphisms in $\mathbb{KOM}(\mathbf{Ban})$
$$\begin{CD}
\dots @>>> 0@>>> 0@>>> X @>>> 0 @>>> 0 @>>> \dots\\
&&@VVV @VVV @VVV @VVV @VVV\\
\dots @>>> 0@>>> Y@>>> 0 @>>> 0 @>>> 0 @>>> \dots
\end{CD}
$$
or, in the language of Section 15, to roofs $\mathcal Q(\text{ morphism\; of\; complexes})\mathcal Q(\text{quasi-isomorphism})^{-1}$
$$\begin{CD}
\dots @>>> 0@>>> 0@>>> X @>>> 0 @>>>  \dots\\
&&@AAA @AAA @AAA @AAA  \text{quasi-isomorphism}\\
\dots @>>> 0@>>> Y@>>> Z @>>> 0 @>>>  \dots\\
&&@VVV @VVV @VVV @VVV \text{ morphism\; of\; complexes}\\
\dots @>>> 0@>>> Y@>>> 0 @>>> 0 @>>>  \dots
\end{CD}
$$

In fact, carrying away this idea, one could form a new category $\mathbf D(\mathbf {Ban})$, whose objects are pairs $(X,n)$ with $X$ a Banach space and $n$ a natural number. The space $\Hom ((X,n), (Y,m))$ is declared to be (the usual)
$\Ext^{n-m}(X,Y)$. This yields a larger category into which $\mathbf {Ban}$ can be embedded in the form $X\to (X,0)$ since also
$\mathfrak L(X,Y)= \Hom ((X,0), (Y,0))$. Of course, $\Ext_{\mathbf {Ban}}^n(X,Y) = \Hom_{\mathbf D(\mathbf {Ban})}(X[0], Y[n])$. Be as it may, from this we infer:

\begin{proposition} $\mathfrak  L_X^{(n)} (Y) = \Ext_{\mathbf {Ban}}^n(X,Y)$
\end{proposition}

According to what is currently known (see \cite[1.8.1, 4.5]{hmbst}) a semi Banach structure on $\Ext_{\mathbf {Ban}}(X,Y)$ is the best one can get. The advantage of direct approaches, such as those in \cite{hmbst,constants,cabecastlong}, is that they provide explicit forms for the seminorm in $\Ext$. However, they do not however provide insights about what occurs with $\Ext^n$. What we did in this paper shows that $\Ext^n_{\mathbf {Ban}}(X,Y)$ carries a semi-Banach topology too (semi-$p$-Banach or semi-quasi-Banach in the more general context of $p$-Banach and quasi Banach spaces).

\section{The heart of Banach spaces}\label{heart}

A category is called additive when the sets $\Hom(A,B)$ have an Abelian group structure compatible with composition, in the sense that $f(g+g')h=fgh + fg'h$, that there is a $0$ object that finite products and coproducts exist. An additive category is Abelian when every arrow has kernel and cokernel and, moreover, each monic arrow is the kernel of its cokernel and every epic arrow is the cokernel of its kernel. Banach and quasi Banach spaces are additive categories but not Abelian and the reason is that the cokernel of an operator $\tau: A\to B$ is the quotient map $B\to B/\overline{\tau[A]}$ (the proof is easy and can be seen in  HMBST, Section 2.4). But that is not the cokernel in the underlying category $\mathcal V$ of vector spaces! In Banach spaces monic arrows are the injective operators and epic arrows are the dense range operators, and that is the reason why even $\mathbb R$ is not injective or projective, and therefore injective or projective objects do not exist (all this was observed by Pothoven \cite{potho}). In particular, it is impossible that every monic (injective) operator is the kernel of its cokernel, as any injective operator with dense range shows.

Still, the reason is so subtle that categorical and homological constructions in $\mathbf {Ban}$ and $\mathbf {QBan}$  mostly work, one way or another ($\mathbb R$, after all, does behave \emph{like} an injective and projective space, doesn't it?). To put a name to his phenomenon Quillen \cite{quill2} introduces the notion of \emph{exact category}: an additive category in which one determines a class of ``short exact sequences" satisfying some suitable conditions \cite[1.1.4]{bbd}, that he calls admissible. There are other notions of exact category: Barr \cite{barr} introduces a notion of exactness so that $$\mathrm{(Abelian) = (exact) + (additive)}$$
Thus, Barr's notion is intrinsic, while Quillen's is extrinsic, in the sense that, \emph{first}, one is already working in an additive category and \emph{then} specifies which are the exact sequences. A detailed exposition of Quillen's notion can be enjoyed in \cite{buhler}. So, Banach and quasi Banach spaces are exact categories in Quillen's sense; why? Because the items to which we award the certification of short exact sequences are ... well ... \emph{the} short exact sequences. Such is the magic that the open mapping theorem brings in $\mathbf {Ban}$ and $\mathbf {QBan}$. Of course that the work of categorical people in defining the notions is much more delicate (and also that of functional analysis people, see \cite{frerik}), because in \emph{a} category no such magic exist (even if the  $\Hom(A,B)$ groups are made of morphisms). In Quillen's definition, admissible exact sequences are made with admissible monics (injective maps) and admissible (epics) subject to some rules.  Just to give you something to think about, try to specify which are the admissible monic and epic in normed spaces, and which are the admissible short exact sequences. Exactly: those.

A detailed list of the rules to follow can be seen in B\"uhler \cite[Definition 2.1]{buhler}, and it reads like this: admissible arrows form a closed class under composition, of course that identities are admissible, and pushouting (pullbacking) admissible monics (epics) yields admissible monics (epics). Things become clearer and nicer if we formulate them as in B\"uhler \cite[Proposition 2.19]{buhler}

\begin{proposition} Consider a commutative square
$$\xymatrix{
A\ar[r]\ar[d] & B\ar[d]\\
C\ar[r]& D}$$
in which the horizontal arrows are admissible monics. The following assertions are equivalent:
\begin{enumerate}
\item The diagram is a pushout.
\item The diagonal sequence
$$\begin{CD} 0@>>> A @>>> C\oplus B @>>> D@>>> 0\end{CD}$$
is exact.
\item The square is a Dolittle diagram: i.e., it is both a pushout and a pull-back diagram.
\item The square is part of a commutative diagram
$$\xymatrix{
0\ar[r]&A\ar[r]\ar[d] & B\ar[d]\ar[r]& E\ar@{=}[d]\ar[r]&0\\
0\ar[r]& C\ar[r]& D\ar[r]& E\ar[r]&0}$$
with exact rows.
\end{enumerate}
\end{proposition}

What we feel when looking at those conditions in $\mathbf{Ban}$ or $\mathbf{QBan}$ is: ok, exactly \emph{that} is what happens! Well, that is why $\mathbf {Ban}$ is an exact category.\\

There have been more efforts to formalize what makes $\mathbf {Ban}$ good enough to be good and correct what makes it bad enough to be bad: that only some maps have ``good" cokernel. The basic idea is: how to embed $\mathbf {Ban}$ as a full subcategory of an Abelian category? Paraphrasing \cite{bbd}:

\adef A heart of a category is an Abelian category that contains it as a full subcategory.\zdef

Thus, guided by Tom Waits let us look for the heart of Banach spaces. The name of Waelbroeck should have appeared first in this regard. His attempt was the construction of the Abelian category of \emph{quotient Banach spaces} \cite{waelb1,waelb2,waelb4,waelb6} and his reasons to do so are described by himself in \cite{waelb3,waelb5}. I find a few difficulties to describe Waelbroeck's work, one of then that it does not help that definitions change from one paper to another. The final Section 4 of \cite{heart} contains a respectful discussion on this topic. Let me focus on the description provided in the paper \cite{waelb1} instead of the much more complete compendium \cite{waelb6}. Waelbroeck's definition in \cite[Def. 1]{waelb1} is ``A quotient Banach space $E/F$ is the quotient of a Banach space $(E, \|\cdot\|_E)$ by a vector subspace $F\subset E$ that is a Banach space with a norm $(F, \|\cdot\|_F)$ finer than the restriction of $\|\cdot\|_E$. In this definition, as far as I understand it, the vector space $E/F$ is part of the definition. This is probably due to a kind of atavic pulsion to ensure that objects of the new category are sets. Once this drive passes and one accepts that the objects of the new category are not sets, Waelbroeck \cite[Def.2.1]{wdef} sets (see also, in \cite[Section 2.1]{waelb6}):

\adef A $q$-space $X|X'$ is a couple $(X, X')$ where $X$ is a Banach space and $X'$ is a vector subspace of $X$ that is a Banach space with a norm $\|\cdot\|_{X'}$ finer than $\|\cdot\|_X$.\zdef

Acting this way things are ok, and we are on our way to defining Wegner's \textbf{Mon} category \cite{heart}. A morphism $\overline {u}: X|X' \To Y|Y'$ between two quotient Banach spaces is in \cite{waelb1} a linear map $\mathfrak u:X/X' \To Y/Y'$ for which there is a bounded operator $u: X\to Y$ fitting in a commutative diagram
\begin{equation}\label{strictd}\xymatrix{
X'\ar[r]\ar[d]_{u'}& X \ar[r]^{Q_1}\ar[d]^u&X/X'\ar[d]^{\mathfrak u} \\
Y'\ar[r]& Y \ar[r]^{Q_2}& Y/Y'} \end{equation}
These morphisms have been called \emph{strict operators}. In \cite[Def. 2.1.2]{waelb6} the author considers that the right equality notion has not been defined. This is important because we need to know what does exactly mean $\overline u=0$: that $\mathfrak u=0$? That $u=0$? And it is not explicitly defined probably because it is not easy to do, in those terms. Indeed, what one is is desperately attempting to say is that homotopy is the right equivalence relation to be set. Homotopy involving linear maps is, however, a sticky affair. Later in \cite{waelb6}, where a quotient Banach space is just an injective operator $f: X'\To X$ between Banach spaces, a strict morphism $f \To g$ is just a couple of operators
\begin{equation*}\label{strict}\xymatrix{
 X'\ar[r]^{u'}\ar[d]_{f}& Y' \ar[d]^g\\
X\ar[r]_u& Y } \end{equation*}
and $(u,u')=0$ if and only if there is an operator $h: X\To Y'$ making the diagram
\begin{equation*}\label{strict}\xymatrix{
 X'\ar[r]^{u'}\ar[d]_{f}& Y' \ar[d]^g\\
X\ar[r]_u\ar[ur]_h& Y } \end{equation*}
commutative. It is clear moreover that $(u,u')=(v,v')$ if and only if $(u-v, u'-v')=0$. Waelbroeck goes on and defines in \cite{waelb1} a \emph{pseudo-isomorphism} when the linear map $b$ is so that the operator $u$ in the diagram (\ref{strictd}) is surjective and $u^{-1}[Y']=X'$
(to be precise, assuming that $f,g$ are as before,  $u^{-1}(g[Y'])=f[X']$). Let's us continue being imprecise for the sake of clarity. A pseudo-isomorphism is necessarily bijective: $0= \mathfrak u(x+X')= ux + Y' \Rightarrow ux\in Y' \Rightarrow x \in X' \Rightarrow x+ X'= y + X' = 0$. The category $\mathbf {qB}$ in \cite{waelb1} has quotient Banach spaces as objects and as morphisms compositions $s^{-1}b$ of strict morphisms $\overline{b}$ and inverses of pseudo-isomorphisms $\overline{s}$. In other words, it is the localization of (quotient Banach spaces, strict morphisms) with respect to pseudo-isomorphisms. The localization is needed here because strict
morphisms do not suffice to get an Abelian category. As we have already mentioned, [71] clearly explains that $\mathbf{qB}$ is one of those categories whose objects are not sets and have no points and its morphisms are not maps. The crucial arrival point is \cite[Lemma 5]{waelb1}:

\begin{lemma}\label{kerqB} Every morphism in $\mathbf {qB}$ has kernel and cokernel\end{lemma}
\begin{proof} As always, the cokernel is the tricky part. Waelbroeck hints the proof like: use that if $\overline{u}: X|X' \To Y|Y'$ is a strict morphism then there exists a cokernel $\mathfrak c: Y/Y'\to (Y/Y')/b[X/X']$ for $\mathfrak u$ in \textbf{Vect}. Back to diagram (\ref{strictd}), $\mathfrak u[X/X']= \mathfrak u Q_1[X]=Q_2u[X]= (u[X]+Y')/Y'$ so that $(Y/Y')/\mathfrak u[X/X'] = (Y/Y')/(u[X]+Y'/Y') = Y/u[X]+Y'$. The space $Y' + u[X]$ can be provided with a Banach structure given by the norm $\|z\| = \inf \{ \|y'\|_{Y'} + \|x\|_X\}$ where the infimum is taken over all representations
$z = y' +  u(x)$. Thus, one can form the object $Y|( Y' + u[X])$ with associated strict morphism $\overline{c}: Y|Y' \to Y| (Y' + u[X])$ and form the commutative diagram
\begin{equation}\label{coke}\begin{CD}
Y'@>>> Y @>>> Y/Y'\\
@VVV @| @VV{\mathfrak c}V\\
Y'+u[X]@>>> Y @>>>Y/ ( Y'+u[X])\end{CD}\end{equation}
Obviously $\mathfrak c \mathfrak u=0$ and therefore $\overline{c}\overline{u} =0$; and if $\overline{w}: Y|Y'  \To S|S'$ is another strict morphism
$$\begin{CD}
Y'@>>> Y @>Q_2>> Y/Y'\\
@VV{w'}V @VVwV @VV{\mathfrak w}V\\
S'@>>> S@>q>>S/S'\end{CD}$$
such that $\overline{w}\overline{u}=0$ then $\mathfrak w$ factors through $\mathfrak c$ in the form $\mathfrak w= \mathfrak m \mathfrak c$. Look at the diagram
$$\xymatrix{
X'\ar[r]\ar[d]& X\ar[d]_u\ar[r]^{Q_1}& X/X'\ar[d]^{\mathfrak b}\\
Y'\ar[r]\ar[d]\ar[ddr]& Y \ar[r]^{Q_2}\ar@{=}[d]\ar[ddr]_{w}& Y/Y'\ar[d]^{\mathfrak c}\ar[ddr]^{\mathfrak w}\\
Y'+u[X]\ar[r]\ar[dr]& Y \ar[r]^-{q_2}\ar[dr]_{m}& Y/  Y'+u[X]\ar[dr]_{\mathfrak m}\\
&S'\ar[r]& S\ar[r]_q& S/S'}$$
Set $m=w$ to get $qm=qw=\mathfrak w Q_2= \mathfrak m\mathfrak cQ_2=\mathfrak mq_2$.
Making the general case $s^{-1}b$ in which $b$ is a strict morphism and $s$ a pseudo-isomorphism is now merry joy: its cokernel is
$cs$.\end{proof}

By construction it takes now little effort \cite[Proposition 6]{waelb1} to show that monic arrows are the kernels of their cokernels and epic arrows are cokernels of their kernels, which leads us to the inexorable conclusion:

\begin{proposition} The category $\mathbf {qB}$ is Abelian.
\end{proposition}

The natural inclusion of $\mathbf {Ban}$ into $\mathbf{qB}$ is to regard $X$ as $X|0$ and thus $\mathbf {Ban}$ is a full subcategory of an Abelian category. The category $\mathbf {qB}$ seems to be a kind of portable derived category of $\mathbf {Ban}$ since, as we said before, it is a localization of the category with objects $F|E$ and strict morphisms with respect to pseudo-isomorphism. In conclusion, Waelbroeck's category $\mathbf {qB}$ is
a heart of Banach spaces. Wegner \cite{heart} cuts the problem straight to arrive to the heart of Banach spaces: he defines the category $\mathbf{Mon}$ whose objects are injective operators $f: X'\to X$ between Banach spaces while an arrow $(\alpha', \alpha): f\to g$ between two objects is a commutative diagram
$$\xymatrix{
X'\ar[r]^{\alpha'}\ar[d]_f& Y'\ar[d]^g\\
X\ar[r]_\alpha& Y}$$
with the following equality notion between morphisms: $(\alpha', \alpha) =0$ if and only if there exists an operator $h$ making a commutative diagram
$$\xymatrix{
X'\ar[r]^{\alpha'}\ar[d]_f& Y'\ar[d]^g\\
X\ar[r]_\alpha \ar[ur]^h& Y}$$

Since $f$ and $g$ are injective it is enough to ask $gh=\alpha$ to get $hf=\alpha'$ free of charge. And therefore, $(\alpha', \alpha) =  (\beta', \beta)$ if and only if $(f'-g', f-g)=0$. If the process that brought us thus far has not vanished into oblivion the reader will realize that
$(\alpha', \alpha) =  (\beta', \beta)$ if and only if $(\alpha', \alpha)$ and $(\beta', \beta)$, seen as morphisms of complexes are homotopic.
The proof only requires to tilt your head (or rotate  $90^\circ$ the page), this way:
$$\xymatrix{\ar[r]&0\ar[r]& X'\ar[r]^{f}\ar[d]_{\alpha'}\ar[dl]& X\ar[d]^g\ar[dl]_h\ar[r]&0\ar[r]\ar[dl]& \\
\ar[r]&0\ar[r]& Y'\ar[r]_{\alpha} & Y\ar[r]&0\ar[r]&}$$

The same rotation shows that Wegner's category is Waelbroeck's category $\mathbf{qB}$ prior to localization. The difference, or so it seems, between Waelbroeck's localized category \textbf{qB} and Wegner's localization of \textbf{hMon[D]} are that i) Waelbroeck does not clearly sets the localizing class; ii) even if he does, it is not clear that the definition of such class is correct. Wegner localizes using a larger system (pulations, see below) than Waelbroeck (pulations whose induced linear application on the quotients is surjective (?)). Indeed, while Wegner shows that if something equivalent to a pulation is a pulation there remains the question of whether a strict morphism equivalent to a pseudo-isomorphism a pseudo-isomorphism. Keeping all that in mind, one of course has:

\begin{lemma} Every morphism in $\mathbf{Mon}$ has kernel and cokernel\end{lemma}
\begin{proof} The cokernel stuff is crystal clear now: given an arrow $(\alpha', \alpha): f\To g$ form the object
$\imath: g[Y'] + \alpha[X] \To Y$ where $g[Y'] + \alpha[X] $ is endowed with the complete norm $\|z\| =\inf \{\|y'\|_{Y'} + \|x\|_X : z = g(y') + \alpha(x)\}$ and the arrow $(g, 1): g\To \imath$
$$\xymatrix{
X'\ar[r]^{\alpha'}\ar[d]_f& Y'\ar[d]^g\ar[r]^-{g}& g[Y'] + \alpha[X] \ar[d]^\imath\\
X\ar[r]_\alpha & Y\ar@{=}[r]_1&Y}$$
One has $(g\alpha', 1\alpha)=0$ since there is the operator $\alpha: X\To g[X] + \alpha[X] $.\end{proof}
Astute readers and Who fans won't get fooled again: this is Lemma \ref{kerqB} verbatim. Now, to proceed with localization, Wegner makes a clever move, full of fun for us. To understand why it would be good to make a detour to review the Kaijser-Pelletier approach to interpolation (see a passing mention in HHI, p.118-119) as exposed in \cite{kaij1,kaij2}. The core idea in those papers is that of Doolittle diagram; namely, a commutative diagram
$$\begin{CD}
A@>>> B\\
@VVV @VVV\\
C@>>> D\end{CD}$$
that is simultaneously a pullback and a pushout diagram. Wegner adopts in \cite[Definition 5]{heart} the name ``pulation" for such square,
introduced by J. Ad\'amek, H. Herrlich, and G. E. Strecker in [Abstract and concrete categories: the joy of cats, Repr. Theory
Appl. Categ. (2006), no. 17, 1--507] p. 205. Even if \emph{the joy of cats} is a wonderful subtitle for a book on categories, pulation is an  awful name at different levels and in different languages. The author, Wegner, is aware that ``other naming conventions (Doolittle diagram, push-me pull-you diagram or bicartesian square) are mentioned in the literature."\\

Let's leave this sad episode in the past and move ahead. Wegner proceeds to localize with respects to the class $D$ of maps $(\alpha',\alpha)$ that provide Doolittle diagrams. After that, Wegner forms the localized category $\heartsuit = \mathbf{Mon}[D]$ and proves the inevitable result that $\heartsuit$ is an Abelian category in which $\mathbf {Ban}$ embeds as a full subcategory in the obvious way: making $X$ correspond to $[0\to X]$; namely:

\begin{proposition} $\heartsuit$ is a heart of Banach spaces.
\end{proposition}

Indeed, categories, like people, may have many hearts. A nice one, intended to provide a solid foundation for Waelbroeck's category of $q$-Banach spaces, was presented by No\"el \cite{noe}. We cannot resits to sketch it with its many virtues:

\subsection*{No\"el's category}  As it is mentioned in \cite[2.1.2]{waelb6} ``The Abelian category of quotient spaces was defined by Waelbroeck in 1962, but it was difficult to work with it. With the ``Miracle Functor", No\"el [64], defined a
similar Abelian category in 1969." What is this miracle functor? Well, as clear as ever, \cite[2.3.3]{waelb6} sets something between
$S \mapsto \mathfrak L(\ell_1(S), X)$ \cite[2.3.6]{waelb6} when $X$ is a Banach space and $X|Y \mapsto \mathfrak L(\ell_1(S), X)/\mathfrak L(\ell_1(S), Y)$ \cite[2.3.11]{waelb6}. See also below. Let $\mathbf K$ be the category whose objects are the Banach spaces $\ell_1(S)$ with $S$ any set (yes, there is some set-theoretic problem here that can somehow be circumvented \cite{noe}).
A $\mathbf Q$-space is an additive contravariant functor $\mathcal F: \mathbf K \To \mathbf{Vect}$. A morphism between two
$\mathbf Q$-spaces is a natural transformation between the functors, and the category obtained this way is called $\mathbf{QESP}$.

\begin{theorem} $\mathbf{QESP}$ is a heart of Banach spaces.\end{theorem}

This means that $\mathbf{QESP}$ is an abelian category and that there is a faithful representation functor $\delta^\bullet: \mathbf {Ban} \To \mathbf{QESP}$, $X \rightarrowtail \mathfrak L^X$ (we already mentioned in Section \ref{catCD} that this functor $\mathbf{Ban}\To \mathbf{Ban}^\mathbf{Ban}$ was a faithful representation). As we have already mentioned, objects of $\mathbf{qB}$ are not sets with elements. No\"el explicitly mentions this in his final remarks (translated from French): ``observe that a $\mathbf Q$-space is not a vector space endowed with an additional structure". However, No\"el observes that there is a natural functor $\sigma: \mathbf{QESP}\To \mathbf{Vect}$ given by
$\mathcal F  \rightarrowtail  \mathcal F(\mathbb C)$ (No\"el works with complex Banach spaces), which makes perfect sense: $\sigma(\mathfrak L^X)= \mathfrak L^X(\mathbb C) = \mathfrak L(\mathbb C, X)= X$. The hardcore of No\"el ideas is \cite[Theorem 2]{noe}:
\begin{proposition} Every $\mathbf Q$-space $\mathcal F$ fits in an exact sequence $0\To \mathfrak L^Y \To \mathfrak L^X \To \mathcal F \To 0$
in which both $Y,X$ are Banach spaces.\end{proposition}whose proof is not, unfortunately, in \cite{noe}. The proposition yields the identity
$\mathcal FX = \mathfrak L(\ell_1(S), X)/\mathfrak L(\ell_1(S), Y)$ (the miracle functor) in $\mathbf{Vect}$ and, in particular,
$\mathcal F(\mathbb C) = X)/Y $. Once one believes that, the rest is easy. Moreover, No\"el claims that $\delta^\bullet$ admits a left adjoint $\mathcal N \xdash \delta^\bullet$, which means that there exists a covariant functor $\mathcal N:  \mathbf{QESP} \To \mathbf {Ban}$ such that
$$\mathfrak L( \mathcal N(\mathcal F), X) = [ \mathcal F, \delta^\bullet(X)] = [\mathcal F, \mathfrak L^X].$$
Since $\mathcal N(\mathcal F)$ has to be a Banach (not a mere vector) space, No\"el claims: $\mathcal N(\mathcal F) = X/\overline Y$.

\subsection*{About the construction of Clausen and Scholze} This issue, that objects of a heart of Banach spaces are not sets remains in the construction of Clausen and Scholze we mention next. Peter Scholze mentioned to us in a private conversation that ``embedding Banach spaces into an Abelian category with wonderful categorical properties has been something that Clausen and I have spent a lot(!) of time thinking about, and we believe we have now
found an excellent way. It is explained in my lectures notes on Analytic
Geometry, www.math.uni-bonn.de/people/scholze/Analytic.pdf . We build a
variant of "complete locally convex topological vector space" that is an
Abelian category (so taking cokernels is fine, even if the image is not
closed), closed under extensions (so we actually have to weaken local convexity, because of the Ribe extension etc. -- as you say, one feels
the need to enlarge Banach spaces to quasi-Banachs), and with further really nice properties (e.g., with compact projective generators). To
prove that it works, we have to prove a more general statement of "arithmetic" flavour -- to set up functional analysis over the reals, we
have to work with the integers! See also  https://xenaproject.wordpress.com/liquid-tensor-experiment/ and
https://mathoverflow.net/questions/386796/nonconv".\\

We end this section with a quote from \cite{waelb6}: ``Mathematicians want to put their hands on elements of an object".
However, those hearts of Banach spaces have no elements. In contrast, S. Lubkin's beautiful Exact Embedding Theorem \cite{lubk} (see also \cite[1.6]{kato}) says:

\begin{theorem} Every Abelian category (whose objects form a
set) admits an additive imbedding into the category of Abelian groups which
carries exact sequences into exact sequences.\end{theorem}

\section{Derivation in Banach spaces, what if....?}

We are ready to face the good, the bad and the ugly issue about derivation in Banach spaces: the same reason why \textbf{Ban} is not Abelian is what makes derivation of $\mathfrak L$ functors work ... as functors $\mathbf{Ban}\to \mathbf{ Vect}$ and not as functor  $\mathbf{Ban}\to \mathbf{Ban}$. Let us explain this assertion: the definition of the homology group at step $n$ of a complex is clear

$$\ker d_{n+1}/ \mathrm{im} \; d_n$$

Here there is something we overlooked in all our previous discussions: when the functor one is trying to derive takes values in a reasonable
category whose objects are sets with some additional structure and the arrows are functions, the meaning of $\mathrm{im} \; d_n$ is clear. So, while deriving $\mathfrak L$ functors in Banach spaces, we are still on solid ground. But, at the cost of having to deal with the image of an operator, something that provides the ``wrong" cokernel and being thus responsible for not having homology groups in the right category (\textbf{Ban}). As we have already mentioned:

\adef\label{image} In an Abelian category, the image of an arrow is the kernel of its cokernel.\zdef

So things work again in an Abelian category and, moreover, Homology ``groups" are objects of the category so that both the construction of the derived category and the final particularizations (1) and (2) after Definition 10 of the derived functor to a functor between the original categories make sense. And, in general? This is a question for categorical people. In the search of a manageable definition for ``image of an arrow" Wegner \cite{heart} isolates the notion of \emph{range}.

\adef The range of a morphisms $f : X \to Y$ is monomorphism $r : R \to Y$ for which there is a morphism $q : X \to R$ in such a way that
 $f = r  q$ and such that for all monomorphism $s: S \to Z$  and morphisms $g : Y \to Z, h: X \to S$ such that $g  f = s  h$ there exists a unique $g' : R \to S$ such that $g  r = s  g'$.\zdef

\noindent which is obviously adapted so that an operator $T:X\To Y$ has range $\overline{T[X]}$. Then, Wegner acknowledges the work of Waelbroeck  formulating:

\adef An additive category $\CAT A$ is a Waelbroeck category if there
exists an additive functor $F : \CAT A \To \CAT {Ab}$ satisfying  the following three conditions:
\begin{enumerate}
\item[(W1)] Every morphism has a kernel and $F$ preserves
kernels
\item[(W2)] Every morphism has a range and $F$ preserves
ranges.
\item[(W3)] The functor $F$ preserves and reflects kernel-cokernel pairs, i.e., $(f, g)$ is a kernel-cokernel
pair in $\CAT A$ if and only if $(Ff,Fg)$ is a kernel-cokernel pair in $\CAT{Ab}$.
\end{enumerate}\zdef

to finally arrive to \cite[Proposition 13]{heart}:

\begin{proposition} $\mathbf{Ban}$ is a Waelbroeck category.
\end{proposition}

that can be completed to\\

\textbf{Coda.} \emph{$\mathbf{QBan}$ is a Waelbroeck category}.

\section{Where are we now?}

Two questions have been wheeling around this survey: How to derive a Banach functor? How to derive a quasi Banach functor? With the elements displayed in this paper there are four possibilities to deal with the first of those questions:

\begin{itemize}
\item Accept that derivation is as it is, including that the derived functor of $\mathfrak L$ is not a Banach functor. There is however a reasonable consolation that they are semi Banach functors.
\item Carry the previous idea to its natural conclusion and derive the semi Banach functor $\mathfrak L: \mathbf{sBan}\to \mathbf{sBan}$.
\item Carry the derivation process strictly (the cokernel of an operator $T:Y\To X$ in $\CAT {Ban}$ is $X/\overline{T[Y]}$) and be prepared to whatever new $\Ext$ functors may appears since the first cohomology group of $X$ at $E$ would be now $\mathfrak L(E, \ell_\infty/X)/\overline{Q \mathfrak L(E, \ell_\infty)}$ where $0\To X\To \ell_\infty \stackrel{Q}\To \ell_\infty/X\To 0$ is a projective presentation of $X$. Now live with that.
\item (suggested by F\'elix Cabello) Embed $\CAT{Ban}$ as a full subcategory of an Abelian category and derive there. Namely, derive in a heart category for Banach spaces. Quoting again Waelbroeck \cite{waelb6}: ``It is not difficult to define categories of quotient spaces and prove that they are Abelian. It is more difficult to develop Functional Analysis in them".
\end{itemize}

A fifth and a six possibilities, if they exist, will hopefully be developed elsewhere.

\end{document}